\documentclass[12pt,reqno]{amsart}

\usepackage{amsmath, amssymb, amsfonts, amsbsy, amsthm, latexsym, graphicx}

\theoremstyle{definition}
\newtheorem{theorem}{Theorem} [section]
\newtheorem{corollary}[theorem]{Corollary}
\newtheorem{lemma}[theorem]{Lemma}
\newtheorem{proposition}[theorem]{Proposition}
\newtheorem{definition}[theorem]{Definition}

\newtheorem{remark}[theorem]{Remark}
\newtheorem{example}[theorem]{Example}
\numberwithin{equation}{section}

\newcommand{\EQ}{\; = \;}
\newcommand{\GE}{\; \ge \;}
\newcommand{\GT}{\; > \;}

\newcommand{\LE}{\; \le \;}
\newcommand{\LT}{\; < \;}
\newcommand{\SUBSET}{\; \subset \;}
\newcommand{\plus}{\; + \;}
\newcommand{\minus}{\; - \;}

\newcommand{\BD}{D_B}

\newcommand{\bigbracket}[1]{\bigl[#1\bigr]}

\newcommand{\C}{\mathbf{C}}
\newcommand{\Cc}{{\mathcal{C}}}
\newcommand{\CHI}{\hbox{\raise .4ex \hbox{$\chi$}}}
\newcommand{\clspan}{{\overline{\mbox{\rm span}}}}

\newcommand{\dist}{{\mathrm{dist}}}
\newcommand{\Ec}{{\mathcal{E}}}
\newcommand{\tEc}{{\tilde{\mathcal{E}}}}
\newcommand{\te}{{\tilde{e}}}
\newcommand{\eps}{\varepsilon}
\newcommand{\Fc}{{\mathcal{F}}}
\newcommand{\tFc}{{\tilde{\mathcal{F}}}}
\newcommand{\tf}{{\tilde{f}}}
\newcommand{\Frac}{{\mathrm{Frac}}}
\newcommand{\Gc}{{\mathcal{G}}}
\newcommand{\tGc}{{\tilde{\mathcal{G}}}}
\newcommand{\tg}{{\tilde{g}}}
\newcommand{\ip}[2]{\langle#1,#2\rangle}

\newcommand{\Bigip}[2]{\Bigl\langle #1, \, #2 \Bigr\rangle}

\newcommand{\Int}{{\mathrm{Int}}}
\newcommand{\cM}{{\mathcal{M}}}
\newcommand{\one}{\mathbf{1}}
\newcommand{\N}{\mathbf{N}}
\newcommand{\norm}[1]{\|#1\|}
\newcommand{\bignorm}[1]{\bigl\|#1\bigr\|}
\newcommand{\Bignorm}[1]{\Bigl\|#1\Bigr\|}
\newcommand{\biggnorm}[1]{\biggl\|#1\biggr\|}
\newcommand{\bomega}{{\bar{\omega}}}
\newcommand{\tomega}{{\tilde{\omega}}}

\newcommand{\bigparen}[1]{\bigl(#1\bigr)}
\newcommand{\Bigparen}[1]{\Bigl(#1\Bigr)}
\newcommand{\biggparen}[1]{\biggl(#1\biggr)}
\newcommand{\plim}{\operatornamewithlimits{\mbox{$p$}-\mathrm{lim}}}
\newcommand{\Q}{\mathbf{Q}}

\newcommand{\qeddeff}{{\qquad\qed}}
\newcommand{\qlim}{\operatornamewithlimits{\mbox{$q$}-\mathrm{lim}}}
\newcommand{\R}{\mathbf{R}}

\newcommand{\Sc}{{\mathcal{S}}}
\newcommand{\set}[1]{\{#1\}}
\newcommand{\bigset}[1]{\bigl\{#1\bigr\}}
\newcommand{\Bigset}[1]{\Bigl\{#1\Bigr\}}
\newcommand{\Span}{\mathrm{span}}

\newcommand{\bx}{{\bar{x}}}
\newcommand{\tx}{{\tilde{x}}}
\newcommand{\Z}{\mathbf{Z}}

\begin{document}

\title[Density, Overcompleteness, and Localization of Frames, II]
{Density, Overcompleteness, and \\ Localization of Frames. \\
II. Gabor Systems}

\author[R.~Balan, P.~G.~Casazza, C.~Heil, and Z.~Landau]
{Radu~Balan, Peter~G.~Casazza, Christopher~Heil, and Zeph~Landau}

\address{\textrm{(R.~Balan)}
Siemens Corporate Research, 
755 College Road East, 
Princeton, New Jersey 08540 USA}
\email{radu.balan@siemens.com}

\address{\textrm{(P.~G.~Casazza)}
Department of Mathematics,
University of Missouri,
Columbia, Missouri 65211 USA}
\email{pete@math.missouri.edu}

\address{\textrm{(C.~Heil)}
School of Mathematics,
Georgia Institute of Technology,
Atlanta, Georgia 30332 USA}
\email{heil@math.gatech.edu}

\address{\textrm{(Z.~Landau)}
Department of Mathematics R8133,
The City College of New York,
Convent Ave at 138th Street,
New York, New York 10031 USA}
\email{landau@sci.ccny.cuny.edu}

\date{June 14, 2005}



\thanks{
The second author was partially supported by NSF Grants
DMS-0102686 and DMS-0405376.
The third author was partially supported by NSF Grant DMS-0139261.
Some of the results of this paper will be announced, without proofs,
in the research announcement \cite{BCHL05b}.}

\begin{abstract}
This work developes a quantitative framework for
describing the overcompleteness of a large class of frames.
A previous paper introduced notions of localization and approximation between
two frames
$\mathcal{F} = \{f_i\}_{i \in I}$ and
$\mathcal{E} = \{e_j\}_{j \in G}$ ($G$ a discrete abelian group),
relating the decay of the expansion of the elements of~$\mathcal{F}$
in terms of the elements of $\mathcal{E}$ via a map $a \colon I \to G$. 
This paper shows that those abstract results yield an array of new
implications for irregular Gabor frames.
Additionally, various Nyquist density results for Gabor frames are
recovered as special cases, and in the process both their meaning
and implications are clarified.
New results are obtained on the excess and overcompleteness of Gabor frames,
on the relationship between frame bounds and density,
and on the structure of the dual frame of an irregular Gabor frame.
More generally, these results apply both to Gabor frames and to systems
of Gabor molecules, whose elements share only a common envelope of
concentration in the time-frequency plane.

The notions of localization and related approximation properties
are a spectrum of ideas that quantify the degree to which elements
of one frame can be approximated by elements of another frame.
In this paper, a comprehensive examination of the interrelations among
these localization and approximation concepts is made, with most implications
shown to be sharp.
\end{abstract}

\copyrightinfo{}{}

\maketitle

\section{Introduction}

The fundamental structural feature of frames that are not Riesz bases is
the overcompleteness of its elements.
To date, even partial understanding of this overcompleteness has been
restricted to limited examples, such as finite-dimensional frames
or highly structured (``lattice'') frames of windowed exponentials or
of time-frequency shifts (Gabor systems).
Together, the ideas and results presented in this paper and
in \cite{BCHL05a} provide a quantitative framework for describing the
overcompleteness of a large class of frames.
The consequences of these ideas are:
(a)~an array of fundamental new results for frames that hold in a
general setting,
(b)~significant new results for the case of Gabor frames,
as well as a new framing of existing results that clarifies their meaning,
and (c)~the presentation of a novel and fruitful point of view for
future research.

Our approach begins with two frames $\Fc = \set{f_i}_{i \in I}$ and
$\Ec = \set{e_j}_{j \in G}$, where $G$ is a discrete abelian group;
we then introduce a notion of the localization of $\Fc$ with respect to $\Ec$.
The idea of localization is that it describes the decay of the coefficients
of the expansion of elements of~$\Fc$ in terms of the elements of $\Ec$.
To make this notion of decay meaningful, a map $a$ from the index set $I$ into
the index set $G$ is introduced.
With this setup, we establish a remarkable
equality between three seemingly unrelated quantities: certain averages
of $\ip{f_i}{\tf_i}$ and $\ip{e_j}{\te_j}$
of frame elements with corresponding canonical dual frame elements,
which we refer to as \emph{relative measures},
and the density of the set $a(I)$ in $G$ \cite[Thm.~3.4]{BCHL05a}.
This equality between density and relative measure is striking since
the relative measure is a function of the frame elements,
while the density is solely determined by the index set~$I$ and
the mapping $a \colon I \to G$.

Due to the length of this work, it is natural to present it in two parts.
The first part, containing the theoretical and structural results that have
driven the research, appeared in \cite{BCHL05a}
(hereafter referred to as ``Part~I'').
In this paper (the second part) we accomplish the following two
main goals.

\medskip
(1) We apply the theoretical results to the case of Gabor systems
$$\Gc(g,\Lambda)
\EQ \set{M_\omega T_x g}_{(x,\omega) \in \Lambda}
\EQ \set{e^{2\pi i \omega \cdot t} g(t-x)}_{(x,\omega) \in \Lambda},$$
which yields a collection of new results that can be summarized as follows.

\medskip
\begin{enumerate}
\item[\textup{(a)}]
\emph{Functions with time-frequency concentration generate localized
Gabor frames} (Theorem~\ref{M1local}).
We show how the degree of localization of a Gabor frame
is tied to the time-frequency concentration of the generating window function
or ``atom''~$g$.
This alone yields a significant improvement over what was previously known
about the approximation properties of irregular Gabor frames.
We extend these results to more general systems of \emph{Gabor molecules}
whose elements are not required to be simple time-frequency shifts of each
other, but instead need only share a common envelope of concentration about
points in the time-frequency plane.

\medskip
\item[\textup{(b)}]
\emph{Characterization of the dual frame of a Gabor frame}
(Theorem~\ref{dualgaborlocalization}).
We prove that if an irregular Gabor frame is generated by a function $g$
which is sufficiently concentrated in the time-frequency plane
(specifically, $g$ lies in the modulation space $M^1$),
then the elements of the dual frame also lie in $M^1$.
We further prove that the dual frame forms a set of Gabor molecules,
and thus, while it need not form a Gabor frame, the elements do share
a common envelope of concentration in the time-frequency plane.
Moreover, this same result applies if the original frame was only
itself a frame of Gabor molecules.
This greatly extends a recent result of Gr\"ochenig and Leinert \cite{GL04}
which covered only the case of lattice Gabor frames.

\medskip
\item[\textup{(c)}]
\emph{A relationship between density of time-frequency shifts and
inner products of frame elements}
(Theorems~\ref{gabordensity}, \ref{gaborbounds}).
We apply the core abstract result of Part~I,
the Density--Relative Measure Theorem \cite[Thm.~3.4]{BCHL05a}.
This implies a remarkable equality between seemingly unrelated quantities:
the density of the time-frequency shifts of a Gabor frame and certain
averages of inner products between Gabor frame elements and the
canonical dual frame elements.
As a consequence we obtain new relationships between the density of the
index set, the frame bounds, and the norm of the window of the Gabor frame.

\medskip
\item[\textup{(d)}]
\emph{The excess of Gabor frames}
(Theorem~\ref{gaborremoval}).
We show that in any overcomplete Gabor frame or set of Gabor molecules,
a set of frame elements with \emph{positive density} can be removed yet
still leave a frame.

\end{enumerate}

\medskip
(2) We provide a comprehensive examination of the interrelations among the
suite of localization localization and approximation concepts introduced in
Part~I, and in most cases provided counterexamples showing that these
implications are sharp.

\smallskip
\subsection{Outline}
We briefly review some results known for irregular Gabor frames
that are related to the themes of this paper and then discuss our new
contributions.

There is an extensive literature available for
``lattice'' Gabor systems of the form
$\Gc(g,\alpha\Z^d \times \beta\Z^d)$ or
$\Gc(g,A(\Z^d))$, where $A$ is an invertible $d \times d$ matrix.
However, until only recently, very few results were available for
irregular Gabor systems $\Gc(g,\Lambda)$ where $\Lambda$ is an
arbitrary subset of $\R^{2d}$.
Some previous papers on irregular Gabor frames are
\cite{Gro93},
\cite{Lan93},
\cite{RS95},
\cite{Jan98},
\cite{CDH99},
\cite{DH00},
\cite{CFZ01},
\cite{HW01},
\cite{SZ02},
\cite{BCHL03b},
\cite{LW03},
\cite{SZ03},
\cite{Gro04}.
We note that many basic questions remain open even for lattice
Gabor frames.
For example, until \cite{BCHL03a}, it was not known if every overcomplete
lattice Gabor frame $\Gc(g,\alpha\Z^d \times \beta\Z^d)$ that was not a Riesz
basis contained an infinite subset that could be removed yet leave a frame.

In \cite{RS95},
Ramanathan and Steger proved a Nyquist density result for certain
irregular Gabor frames.
Together with the extensions from \cite{CDH99},
this can be stated as follows
(compare also \cite{Lan93}, \cite{GR96}, \cite{Jan98}):

\smallskip
\begin{enumerate}
\item[(a)]
If $\Gc(g,\Lambda)$ is a frame for $L^2(\R^d)$, then
$1 \le \BD^-(\Lambda) \le \BD^+(\Lambda) < \infty$;

\smallskip
\item[(b)]
If $\Gc(g,\Lambda)$ is a Riesz basis for $L^2(\R^d)$,
then $\BD^-(\Lambda) = \BD^+(\Lambda) = 1$.
\end{enumerate}

\smallskip\noindent
Here $\BD^\pm(\Lambda)$ are the lower and upper Beurling densities of
$\Lambda$, which are defined precisely in Example~\ref{gabordensityrel}.
In the case that $\Lambda$ is a rectangular lattice of the form
$\Lambda = \alpha\Z^d \times \beta\Z^d$,
several alternative proofs of these facts are known, see
\cite{Rie81},
\cite{Bag90},
\cite{Dau90},
\cite{Jan94},
\cite{DLL95},
\cite{BR03}.
For additional history and references see the exposition in
\cite{Dau92}, \cite{BHW95}, \cite{Gro01}.
Note that a lattice $\Lambda = \alpha\Z^d \times \beta\Z^d$ has uniform
Beurling densities
$\BD^+(\alpha\Z^d \times \beta\Z^d)
= \BD^-(\alpha\Z^d \times \beta\Z^d)
= 1/(\alpha\beta)^d$.

Ramanathan and Steger's results showed that it is not the algebraic structure
of the lattice as such that is essential, but rather the fact that
Gabor frames satisfy a certain Homogeneous Approximation Property
(essentially the weak HAP presented in this paper, using
as a reference system a Gabor frame with a Gaussian generating function).
For a Gabor frame, the HAP essentially states that the rate of approximation
of a given function by Gabor frame elements is invariant under time-frequency
shifts of the function.
This is remarkable in the case of irregular Gabor frames, since there is no
structure to relate the specific frame elements used in the approximation
of one time-frequency shift $M_q T_p f$ to those in an approximation of
another time-frequency shift of $f$.
The HAP is a fundamental property of Gabor frames, yet very few papers
subsequent to \cite{RS95} have made use of it.

We investigate the localization properties of irregular Gabor frames
in Section~\ref{gaborsec}.
Following the introduction of some notation in Section~\ref{gaborprelims},
we show in Section~\ref{boxlocalsec} that every Gabor system has at least a
weak amount of localization ($\ell^2$-row decay) with respect to the
Gabor orthonormal basis generated by the box function, and as a consequence,
we recover the fact that every Gabor Bessel sequence has finite density.

In Section~\ref{M1localsec} we show that much stronger localization is
obtained if the reference Gabor system is generated by a function that is
well-concentrated in the time-frequency plane.
The degree of concentration is quantified by the modulation spaces, which
are the natural function spaces associated with Gabor analysis.
In particular, we show that if $\Gc(g,\Lambda)$ is an arbitrary Gabor
system such that $g$ lies in the modulation space $M^p$ ($1 \le p \le 2$),
and $\Gc(\phi,\alpha\Z^d \times \beta\Z^d)$ is a reference Gabor system
whose generator $\phi$ lies in the modulation space $M^1$, then
$(\Gc(g,\Lambda), \, a, \, \Gc(\phi,\alpha\Z^d \times \beta\Z^d))$
is $\ell^p$-localized.
Additionally,
$(\Gc(g,\Lambda), \, a)$ is $\ell^1$-self-localized if $p=1$.
In particular, by the implications established in Theorem~\ref{relations},
$\ell^2$-localization for a frame
implies both the Strong HAP and the Strong Dual HAP, so these
results greatly improve the Homogeneous Approximation Property
previously established for Gabor frames in \cite{RS95}, \cite{CDH99}.

In Section~\ref{newresults} we use this extended knowledge of the
localization properties of Gabor systems to derive new implications
for irregular Gabor frames.
In Section~\ref{gabordensitysec} we recover the density results for
Gabor frames given above, and furthermore we extend the meaning of density by
showing that the Beurling densities of $\Lambda$ are related to the measures
of the Gabor frame.
We show that the upper and lower measures of a Gabor frame satisfy
$\cM^\pm(\Gc(g,\Lambda)) = 1/\BD^\mp(\Lambda)$.
This gives a new interpretation of the density of a Gabor frame,
and as a consequence we obtain new relations among the density of the
index set, the frame bounds, and the norm of the generator.
In particular, we show that if $\Gc(g,\Lambda)$ is a tight frame then
the index set $\Lambda$ must have uniform Beurling density, i.e.,
$\BD^-(\Lambda) = \BD^+(\Lambda)$.
Thus tight Gabor frames require a certain amount of uniformity in the
index set.

In Section~\ref{gaborremovalsec},
we show that if $\Gc(g,\Lambda)$ is a Gabor frame whose generator
$g$ lies in $M^1$, then whenever $\BD^-(\Lambda) > 1$, there is
not merely an infinite subset but a subset with positive density which
may be removed from the frame yet still leave a frame.

In Section~\ref{gabordualsec} we address the fundamental question of the
structure of the canonical dual frame of an irregular Gabor frame.
As is well-known, the canonical dual frame of a lattice Gabor frame is again
a lattice Gabor frame, indexed by the same lattice.
This need not be the case for an irregular Gabor frame
$\Gc(g,\Lambda)
= \set{e^{2\pi i \omega \cdot t} g(t-x)}_{(x,\omega) \in \Lambda}$
with an arbitrary index set~$\Lambda$.
A canonical dual frame $\tGc = \set{\tg_{x,\omega}}_{(x,\omega) \in \Lambda}$
will exist, but to date essentially nothing has been known about this
dual beyond the fact that each dual frame element $\tg_{x,\omega}$
is some function in $L^2(\R^d)$.
We prove that if $g$ possesses sufficient time-frequency concentration,
namely that $g \in M^1$, then each $\tg_{x,\omega}$ possesses the same
concentration, i.e., $\tg_{x,\omega} \in M^1$.
For the case of lattice systems ($\Lambda = \alpha\Z^d \times \beta\Z^d$),
this result was previously obtained by Gr\"ochenig and Leinert \cite{GL04}.
However, in addition to extending to the completely irregular setting,
we also prove that the dual frame $\tGc$ is Gabor-like.
Namely, even though $\tGc$ need not itself be a Gabor frame, we show that
it will form a set of \emph{Gabor molecules},
meaning that each function $\tg_{x,\omega}$
is concentrated in the time-frequency plane about the point $(x,\omega)$
with a common envelope of concentration for each $(x,\omega) \in \Lambda$.
In fact, this result holds even if the original frame was only itself
a frame of Gabor molecules.

We conclude the Gabor portion of the paper in Section~\ref{molecules} where we observe that most of the results
obtained for Gabor frames carry over with minor changes to the case of
Gabor molecules.

For simplicity of presentation, most of our results will be stated for the
case of Gabor frames for all of $L^2(\R^d)$, but many of them can be
extended to the case of Gabor frame sequences, or to Gabor frames with
multiple generators, by making use of the machinery developed in Part~I
and this paper.
Results analogous to the ones formulated for Gabor systems could also be
formulated for the case of windowed exponentials.

Finally, in Section~\ref{relationsec} we carefully examine the interrelations between the range of localization properties and approximation concepts that have been used to develop the theory.
 A set of approximation properties for abstract frames introduced in
Part~I is given in Definition~\ref{approxdef}.
These are defined in terms of how well the elements
of the reference system are approximated by finite linear combinations
of frame elements, or vice versa, and provide an abstraction for
general frames of the essential features of the Homogeneous Approximation
Property (HAP) that is known to hold for Gabor frames or windowed exponentials
(see \cite{RS95}, \cite{GR96}, \cite{CDH99}).
 Theorem ~\ref{relations} establishes
an extensive list of implications that hold among
the localization and approximation properties,
and additionally in most cases we provide examples which show that
these implications are sharp.
In particular, there is an equivalence between
$\ell^2$-column decay and the HAP, and between
$\ell^2$-row decay and a dual HAP.

We believe that localization is a powerful and useful new concept.
As evidence of this fact, we note that Gr\"ochenig has independently
introduced a concept of localized frames, for a completely different
purpose, in \cite{Gro04}.
We learned of Gr\"ochenig's results shortly after completion of our own
major results.
The definitions of localizations presented here and in \cite{Gro04}
differ, but the fact that this single concept has independently arisen
for two very distinct applications shows its utility.
In his elegant paper, Gr\"ochenig has shown that frames which are sufficiently
localized in his sense provide frame expansions not only for the Hilbert
space~$H$ but for an entire family of associated Banach function spaces.
Gr\"ochenig further showed that if a frame is sufficiently localized
in his sense (a~polynomial or exponential localization) then the dual frame
is similarly localized.

\smallskip
\subsection{General Notation} \label{prelimsec}

We use the notation from Part~I, which we briefly review here.
$H$~will refer to a separable Hilbert space.
The frame or system of interest will be indexed by a countable index set $I$.
The reference frame or system will be indexed by an additive discrete
group $G$ of the form
$$G \EQ \prod_{i=1}^d a_i \Z \, \times \, \prod_{j=1}^e \Z_{n_j}.$$
We define a metric on $G$ as follows.
If $m_j \in \Z_{n_j}$, set $\delta(m_j) = 0$ if $m_j =0$, otherwise
$\delta(m_j) = 1$.
Then given $g = (a_1 n_1, \dots, a_d n_d, m_1, \dots, m_e) \in G$, set
$$|g|
\EQ \sup\bigset{|a_1n_1|, \, \dots,\, |a_d n_d|, \,
                \delta(m_1), \, \dots, \, \delta(m_e)}.$$
The metric is $d(g,h) = |g-h|$ for $g$, $h \in G$.
The reader can simply take $G = \Z^d$ without much loss of insight
on a first reading (the metric in this case is simply the $\ell^\infty$
metric on $\R^d$ restricted to $\Z^d$).

We implicitly assume that there exists a map $a \colon I \to G$
associated with $I$ and~$G$.
This map will often not be injective.
For each integer $N > 0$ we let
$$S_N(j)
\EQ \Bigset{k \in G : |k-j| \le \frac{N}2}$$
denote a discrete ``cube'' or ``box'' in $G$ centered at $j \in G$.
The cardinality of $S_N(j)$ is independent of $j$.
For example, if $G = \Z^d$ then
$|S_{2N}(j)| = |S_{2N+1}(j)| = (2N+1)^d$.
We let $I_N(j)$ denote the inverse image of $S_N(j)$ under $a$, i.e.,
$$I_N(j) \EQ a^{-1}(S_N(j))
\EQ \set{i \in I : a(i) \in S_N(j)}.$$

\smallskip
\subsection{Notation for Frames and Riesz Bases}

We use standard notations for frames and Riesz bases as found in the texts
\cite{Chr03}, \cite{Dau92}, \cite{Gro01}, \cite{You01}
or the research-tutorials \cite{Cas00}, \cite{HW89}.

A sequence $\Fc = \set{f_i}_{i \in I}$ is a \emph{frame} for $H$ if there
exist \emph{frame bounds} $A$, $B > 0$ such that
$A \, \norm{f}^2 \le \sum_{i \in I} |\ip{f}{f_i}|^2 \le B \, \norm{f}^2$
for all $f \in H$.
The \emph{analysis operator} $T \colon H \to \ell^2(I)$ is
$Tf = \set{\ip{f}{f_i}}_{i \in I}$, and its adjoint
$T^* c = \sum_{i \in I} c_i \, f_i$
is the \emph{synthesis operator}.
The \emph{Gram matrix} is
$T T^* = [\ip{f_i}{f_j}]_{i,j \in I}$.
The \emph{frame operator}
$Sf = T^* T f = \sum_{i \in I} \ip{f}{f_i} \, f_i$
is a bounded, positive, and invertible mapping of $H$ onto itself.
The \emph{canonical dual frame} is
$\tFc = S^{-1}(\Fc) = \set{\tf_i}_{i \in I}$
where $\tf_i = S^{-1} f_i$.
We call $\Fc$ a \emph{tight frame} if we can take $A=B$, and
a \emph{Parseval frame} if we can take $A=B=1$.
If $\Fc$ is any frame, then $S^{-1/2}(\Fc)$ is the
\emph{canonical Parseval frame} associated to $\Fc$.
We call $\Fc$ a \emph{uniform norm frame} if all the frame elements have
identical norms, i.e., if $\norm{f_i} = const.$ for all $i \in I$.
A frame is \emph{norm-bounded below} if $\inf_i \norm{f_i} > 0$.

A sequence which satisfies the upper frame bound estimate,
but not necessarily the lower estimate, is called
a \emph{Bessel sequence} and $B$ is a \emph{Bessel bound}.

A sequence $\Fc = \set{f_i}_{i \in I}$
that is a frame for its closed linear span in $H$
is called a \emph{frame sequence}.
In this case $\tFc = \set{\tf_i}_{i \in I}$ will denote its
canonical dual frame within $\clspan(F)$, and
$P_\Fc$ will denote the orthogonal projection of $H$ onto $\clspan(\Fc)$.

A frame is a basis if and only if it is a Riesz basis.
A \emph{Riesz sequence} is a sequence that forms a
Riesz basis for its closed linear span in~$H$.

\smallskip
\section{Density, Localization, HAP, and Relative Measure} \label{section2}

In this section we recall basic definitions from Part~I,
and show how they are implemented for the case of Gabor systems.
In Theorem~\ref{relations} we will derive an extended set of implications
that hold among these properties, and provide examples in
Section~\ref{relationappend} showing that most of those implications
are sharp.

\smallskip
\subsection{Density} \label{densitysec}

Given an index set $I$ and a map $a \colon I \to G$, we define
the density of $I$ by computing the Beurling density of its image $a(I)$
as a subset of $G$.
Note that we regard $I$ as a sequence, and hence repetitions of images
count in determining the density.

\begin{definition}[Density] \label{densitydef}
The \emph{lower and upper densities of $I$ with respect to $a$} are
\begin{equation} \label{lowerdensity}
D^-(I,a)
\EQ \liminf_{N \to \infty} \inf_{j \in G} \frac{|I_N(j)|}{|S_N(j)|},
\qquad
D^+(I,a)
\EQ \limsup_{N \to \infty} \sup_{j \in G} \frac{|I_N(j)|}{|S_N(j)|},
\end{equation}
respectively.
These quantities could be zero or infinite,
$0 \le D^-(I,a) \le D^+(I,a) \le \infty$.
When $D^-(I,a) = D^+(I,a) = D$ we say $I$ has \emph{uniform density}~$D$.
~\qed
\end{definition}

These lower and upper densities are only the extremes of the
possible densities that we could naturally assign to $I$ with respect to $a$.
In particular, instead of taking the infimum or supremum over all possible
centers in \eqref{lowerdensity} we could choose one
specific sequence of centers, and instead of computing the liminf or limsup
we could consider the limit with respect to some ultrafilter
(some background on ultrafilters is provided in Appendix~\ref{ultraappend}).
The different possible choices of ultrafilters and sequences of centers
gives us the following natural collection of definitions of density.

\begin{definition}
Let $p$ be a free ultrafilter, and let $c = (c_N)_{N \in \N}$ be any
sequence of centers $c_N \in G$.
Then the \emph{density of $I$ with respect to $a$, $p$, and $c$} is
$$D(p,c)
\EQ D(p,c;I,a)
\EQ \plim_{N \in \N} \frac{|I_N(c_N)|}{|S_N(c_N)|}.
\qeddeff$$
\end{definition}

\begin{example}
If $I=G$ and $a$ is the identity map, then $I_N(j) = S_N(j)$ for every $N$
and~$j$, and hence $D(p,c) = D^-(I,a) = D^+(I,a) = 1$
for every choice of free ultrafilter $p$ and sequence of centers~$c$.
\qed
\end{example}

It follows from basic properties of ultrafilters that we always have
$D^-(I,a) \le D(p,c) \le D^+(I,a)$.
It is shown in \cite[Lem.~2.5]{BCHL05a} that
there exist free ultrafilters $p^-$, $p^+$ and sequence of centers
$c^- = (c_N^-)_{N \in \N}$, $c^+ = (c_N^+)_{N \in \N}$ in $G$ such that
$D^-(I,a) = D(p^-,c^-)$
and
$D^+(I,a) = D(p^+,c^+)$.

\begin{example}[Gabor Systems] \label{gabordensityrel}
Consider an arbitrary Gabor system $\Fc = \Gc(g,\Lambda)$
and a reference lattice Gabor system
$\Ec = \Gc(\phi,\alpha\Z^d \times \beta\Z^d)$.
The index set $I = \Lambda$ is a countable sequence of points in $\R^{2d}$,
and the reference group is $G = \alpha\Z^d \times \beta\Z^d$.
A natural map $a \colon \Lambda \to G$,
that we will employ whenever dealing with Gabor systems,
is rounding to a near element of~$G$, i.e.,
\begin{equation} \label{adef1}
a(x,\omega)
\EQ \bigparen{\alpha \, \Int\bigparen{\tfrac{x}{\alpha}},
              \beta \, \Int\bigparen{\tfrac{\omega}{\beta}}},
\qquad (x,\omega) \in \Lambda,
\end{equation}
where $\Int(x) = (\lfloor x_1 \rfloor, \dots, \lfloor x_d \rfloor)$.
With this setup, $S_N(j)$ is the intersection of
$\alpha\Z^d \times \beta\Z^d$ with the cube $Q_N(j)$ in $\R^{2d}$
centered at $j$ with side lengths $N$.
Such a cube contains approximately $N^{2d}/(\alpha\beta)^d$ points of
$\alpha\Z^d \times \beta\Z^d$; precisely,
$$\lim_{N \to \infty} \frac{|S_N(j)|}{N^{2d}} \EQ \frac1{(\alpha\beta)^d}.$$
Likewise, because $a$ is a bounded perturbation of the identity map,
the number of points in $I_N(j)$ is asymptotically the cardinality of
$\Lambda \cap Q_N(j)$.
Consequently, the standard definition of the upper Beurling density
$\BD^+(\Lambda)$ of $\Lambda$ is related to our definition of the
upper density of $\Lambda$ with respect to $a$ as follows:
\begin{align}
\BD^+(\Lambda)
& \EQ \limsup_{N \to \infty} \sup_{j \in \R^{2d}}
      \frac{|\Lambda \cap Q_N(j)|}{N^{2d}} \notag \\[1 \jot]
& \EQ \frac1{(\alpha\beta)^d} \, \limsup_{N \to \infty}
      \sup_{j \in \alpha\Z^d \times \beta\Z^d}
      \frac{|I_N(j)|}{|S_N(j)|}
\EQ \frac1{(\alpha\beta)^d} \, D^+(\Lambda,a). \label{beurlingdef}
\end{align}
Similarly the lower Beurling density of $\Lambda$ is
$\BD^-(\Lambda) = (\alpha\beta)^{-d} \, D^-(\Lambda,a)$.
In particular, when $\alpha\beta = 1$ (the ``critical density'' case),
our definition coincides with Beurling density, but in general the extra
factor of $(\alpha\beta)^d$ must be taken into account.

The map $a$ given in \eqref{adef1} is the one we will use when dealing with
Gabor systems, but any bounded perturbation of $a$ would serve just as well.
That is, given $\delta > 0$, we could map $(x,\omega)$ to any point in
$G = \alpha\Z^d \times \beta\Z^d$ that is within a distance $\delta$ of
$(\alpha \, \Int(\tfrac{x}{\alpha}), \beta \, \Int(\tfrac{\omega}{\beta}))$
without any change in the results.
\qed
\end{example}

\smallskip
\subsection{The Localization Properties} \label{localizationsec}

The words ``column'' and ``row'' in the following definition
refer to the $I \times G$ cross-Grammian matrix
$[\ip{f_i}{e_j}]_{i \in I, j \in G}$.
We think of the elements in locations $(i,a(i))$ as corresponding to the
main diagonal of this matrix.

\begin{definition}[Localization] \label{localizationdef}
Let $\Fc = \set{f_i}_{i \in I}$ and
$\Ec = \set{e_j}_{j \in G}$ be sequences in $H$
and $a \colon I \to G$ an associated map.

\smallskip
\begin{enumerate}
\item[(a)]
We say $\Fc$ is \emph{$\ell^p$-localized} with respect to the
reference sequence $\Ec$ and the map $a$,
or simply that \emph{$(\Fc,a,\Ec)$ is $\ell^p$-localized}, if
$$\sum_{j \in G} \, \sup_{i \in I} \, |\ip{f_i}{e_{j+a(i)}}|^p
\LT \infty.$$
Equivalently, there must exist an $r \in \ell^p(G)$ such that
$$\forall\, i \in I, \quad
\forall\, j \in G, \quad
|\ip{f_i}{e_j}| \LE r_{a(i)-j}.$$

\medskip
\item[(b)]
We say that $(\Fc,a,\Ec)$ has \emph{$\ell^p$-column decay} if for every
$\eps > 0$ there is an integer $N_\eps > 0$ so that
$$\forall\, j \in G, \quad
\sum_{i \in I \setminus I_{N_\eps}(j)} |\ip{f_i}{e_j}|^p \LT \eps.$$

\medskip
\item[(c)]
We say $(\Fc,a,\Ec)$ has \emph{$\ell^p$-row decay} if for every
$\eps > 0$ there is an integer $N_\eps > 0$ so that
$$\forall\, i \in I, \quad
\sum_{j \in G \setminus S_{N_\eps}({a(i)})} |\ip{f_i}{e_j}|^p \LT \eps.
\qeddeff$$
\end{enumerate}
\end{definition}

Note that given a sequence $\Fc$, the definition of localization is dependent 
upon both the choice of reference sequence $\Ec$ and the map $a$.

\begin{example}[Gabor Systems] \label{gaborexample}
For motivation, consider the especially simple case of Gabor systems both
indexed by $\Z^{2d}$, i.e.,
$\Fc = \Gc(g,\Z^{2d}) = \set{M_n T_k g}_{(k,n) \in \Z^{2d}}$
and
$\Ec = \Gc(\phi,\Z^{2d}) = \set{M_n T_k \phi}_{(k,n) \in \Z^{2d}}$.
The map $a \colon \Z^{2d} \to \Z^{2d}$ given by \eqref{adef1} is
the identity map, and
the $\Z^{2d} \times \Z^{2d}$ cross-Grammian matrix
$$\bigbracket{\ip{M_m T_j g}{M_n T_k \phi}}
    _{(j,m) \in \Z^{2d}, \, (k,n) \in \Z^{2d}}
\EQ \bigbracket{\ip{g} {M_{n-m} T_{k-j} \phi}}
    _{(j,m) \in \Z^{2d}, \, (k,n) \in \Z^{2d}}$$
is Toeplitz.
Set $r_{k,n} = |\ip{g}{M_n T_k \phi}|$.
If $\Ec = \Gc(\phi,\Z^{2d})$ is a Bessel sequence with Bessel bound~$B$, then 
$$\sum_{(k,n) \in \Z^{2d}} r_{k,n}^2
\EQ \sum_{(k,n) \in \Z^{2d}} |\ip{g}{M_n T_k \phi}|^2
\LE B \, \norm{g}_2^2
\LT \infty,$$
so $r \in \ell^2(\Z^{2d})$, and hence
$(\Gc(g,\Z^{2d}), \, a, \, \Gc(\phi,\Z^{2d}))$ is $\ell^2$-localized.

Unfortunately, Gabor frames indexed by $\alpha\Z^d \times \beta\Z^d$
with $\alpha\beta=1$ are not very useful in practice.
It can easily be shown via Zak transform techniques that if such a system is
a frame for $L^2(\R^d)$ then it will be a Riesz basis.
However, the \emph{Balian--Low Theorem} (BLT) states that the generator $g$
of such a Gabor Riesz basis cannot be simultaneously well-concentrated in both
time and frequency.
For exposition and references on the BLT see the survey \cite{BHW95}.
Some recent results on the BLT in higher dimensions
are in \cite{BCGP02}, \cite{GHHK02}.

Therefore, consider an arbitrary Gabor system
$\Fc = \Gc(g,\Lambda) = \set{M_\omega T_x g}_{(x,\omega) \in \Lambda}$,
where $\Lambda \subset \R^{2d}$.
For a reference system take a lattice Gabor system of the form
$\Ec = \Gc(\phi, \alpha\Z^d \times \beta\Z^d)
= \set{M_\eta T_u \phi}_{(u,\eta) \in \alpha\Z^d\times\beta\Z^d}$,
where $\alpha$, $\beta > 0$.
We regard $\Ec$ as being indexed by
$G = \alpha\Z^d \times \beta\Z^d$
and use the natural map $a \colon \Lambda \to G$ given in \eqref{adef1}
that sends an element of $\Lambda$ to a near element of~$G$.
It is no longer the case that the cross-Grammian matrix
$\bigbracket{\ip{M_\omega T_x g}{M_\eta T_u \phi}}
 _{(x,\omega) \in \Lambda, \, (u,\eta) \in G}$
is Toeplitz, but we will show in Theorem~\ref{M1local} that if $\Lambda$
has finite density and $\phi$ possesses a certain amount of joint
concentration in time and frequency then
$(\Gc(g,\Lambda), \, a, \, \Gc(\phi, \alpha\Z^d \times \beta\Z^d))$
is $\ell^2$-localized.
The specific requirement on $\phi$ is that it must lie in the modulation
space $M^1(\R^d)$, which is defined precisely in Section~\ref{M1localsec}.
Moreover, the localization can be improved by also imposing a time-frequency
concentration condition on $g$.
Specifically, we show in Theorem~\ref{M1local} that if
$g \in M^p(\R^d)$ and $\phi \in M^1(\R^d)$,
then $(\Gc(g,\Lambda), \, a, \, \Gc(\phi, \alpha\Z^d \times \beta\Z^d))$
is $\ell^p$-localized.
~\qed
\end{example}

\begin{remark}
For comparison, let us give Gr\"ochenig's notion of localization from
\cite{Gro04}.
Let $I$ and $J$ be countable index sets in $\R^d$ that are separated,
i.e., $\inf_{i \ne j \in I} |i-j| > 0$ and similarly for $J$.
Then $\Fc = \set{f_i}_{i \in I}$ is \emph{$s$-polynomially localized}
with respect to a Riesz basis $\Ec = \set{e_j}_{j \in J}$ if
for every $i \in I$ and $j \in J$ we have
$$|\ip{f_i}{e_j}| \LE C \, (1 + |i-j|)^{-s}
\qquad\text{and}\qquad
|\ip{f_i}{\te_j}| \LE C \, (1 + |i-j|)^{-s},$$
where $\set{\te_j}_{j \in J}$ is the dual basis to 
$\set{e_j}_{j \in J}$.
Likewise $\Fc = \set{f_i}_{i \in I}$ is \emph{exponentially localized}
with respect to a Riesz basis $\Ec = \set{e_j}_{j \in J}$ if
for some $\alpha > 0$ we have for every $i \in I$ and $j \in J$ that
$$|\ip{f_i}{e_j}| \LE C \, e^{-\alpha |i-j|}
\qquad\text{and}\qquad
|\ip{f_i}{\te_j}| \LE C \, e^{-\alpha |i-j|}.
\qeddeff$$
\end{remark}

\smallskip
\subsection{The Approximation Properties} \label{approxsec}

The following approximation properties extract the essence of the
Homogeneous Approximation Property that is satisfied by Gabor frames,
but without reference to the exact structure of Gabor frames.

\begin{definition}[Homogeneous Approximation Properties] \label{approxdef}
Let $\Fc = \set{f_i}_{i \in I}$ be a frame for~$H$ with
canonical dual $\tFc = \set{\tf_i}_{i \in I}$, and
let $\Ec = \set{e_j}_{j \in G}$ be a sequence in $H$.
Let $a \colon I \to G$ be an associated map.

\smallskip
\begin{enumerate}
\item[(a)]
We say $(\Fc,a,\Ec)$ has the \emph{weak HAP}
if for every $\eps > 0$, there is an integer $N_\eps > 0$
so that for every $j \in G$ we have
$\dist\bigparen{e_j, \,\,
\clspan\bigset{\tf_i : i \in I_{N_\eps}(j)}}
< \eps$.
Equivalently, there must exist scalars $c_{i,j}$, with only finitely many
nonzero, such that
\begin{equation} \label{weakHAPdef}
\Bignorm{e_j - \sum_{i \in I_{N_\eps}(j)} c_{i,j} \, \tf_i}
\LT \eps.
\end{equation}

\medskip
\item[(b)]
We say $(\Fc,a,\Ec)$ has the \emph{strong HAP}
if for every $\eps > 0$, there is an integer $N_\eps > 0$
so that for every $j \in G$ we have
$$\Bignorm{e_j - \sum_{i \in I_{N_\eps}(j)} \ip{e_j}{f_i} \, \tf_i}
\LT \eps. \qeddeff$$
\end{enumerate}
\end{definition}

\begin{definition}[Dual Homogeneous Approximation Properties]
Let $\Fc = \set{f_i}_{i \in I}$ be a sequence in $H$, and
let $\Ec = \set{e_j}_{j \in G}$ be a frame for $H$ with
canonical dual $\tEc = \set{\te_j}_{j \in G}$.
Let $a \colon I \to G$ be an associated map.

\smallskip
\begin{enumerate}
\item[(a)]
We say $(\Fc,a,\Ec)$ has the \emph{weak dual HAP}
if for every $\eps > 0$, there is an integer $N_\eps > 0$
so that for every $i \in I$ we have
$\dist\bigparen{f_i, \,\,
\clspan\bigset{\te_j : j \in S_{N_\eps}(a(i))}}
< \eps$.

\medskip
\item[(b)]
We say $(\Fc,a,\Ec)$ has the \emph{strong dual HAP}
if for every $\eps > 0$, there is an integer $N_\eps > 0$
so that for every $i \in I$ we have
$\bignorm{f_i - \sum_{j \in S_{N_\eps}(a(i))} \ip{f_i}{e_j} \, \te_j}
< \eps$.
~\qed
\end{enumerate}
\end{definition}

\smallskip
\subsection{Self-Localization} \label{selflocsec}

It is also useful to consider localizations where the system
$\Fc = \set{f_i}_{i \in I}$ is compared to itself or to its canonical
dual frame instead of to a reference system~$\Ec$.
An analogous polynomial or exponential ``intrinsic localization'' was
independently introduced by Gr\"ochenig in \cite{Gro03};
see also \cite{For03}, \cite{GF04}.
Although there is no reference system, we still require a mapping
$a \colon I \to G$ associating $I$ with a group $G$.

\begin{definition}[Self-localization]
Let $\Fc = \set{f_i}_{i \in I}$ be a sequence in $H$, and
let $a \colon I \to G$ be an associated map.

\smallskip
\begin{enumerate}
\item[(a)]
We say that $(\Fc,a)$ is \emph{$\ell^p$-self-localized} if there exists
$r \in \ell^p(G)$ such that
$$\forall\, i, j \in I, \quad
|\ip{f_i}{f_j}| \LE r_{a(i) - a(j)}.$$

\medskip
\item[(b)]
If $\Fc$ is a frame sequence, then we say that $(\Fc,a)$ is
\emph{$\ell^p$-localized with respect to its canonical dual frame sequence}
$\tFc = \set{\tf_i}_{i \in I}$ if there exists $r \in \ell^p(G)$ such that
$$\forall\, i, j \in I, \quad
|\ip{f_i}{\tf_j}| \LE r_{a(i) - a(j)}. \qeddeff$$
\end{enumerate}
\end{definition}

\medskip
\begin{remark} \label{selfremark}
If $I = G$ and $a$ is the identity map, then $(\Fc,a)$ is
$\ell^1$-self-localized if and only if $(\Fc,a,\Fc)$ is $\ell^1$-localized.
However, if $a$ is not the identity map, then this need not be the case.
For example, every orthonormal basis is $\ell^1$-self-localized regardless
of which map~$a$ is chosen, but in Example~\ref{selflocexample} we construct
an orthonormal basis $\Fc = \set{f_i}_{i \in \Z}$
and a map $a \colon \Z \to \Z$ such that $(\Fc,a,\Ec)$ is not
$\ell^1$-localized for any Riesz basis $\Ec$; in fact,
$(\Fc,a,\Ec)$ cannot even possess both $\ell^2$-column decay and
$\ell^2$-row decay for any Riesz basis $\Ec$.
\qed
\end{remark}

We show in Example~\ref{selfexample} that $\ell^1$-localization with 
respect to the dual frame does not imply $\ell^1$-self-localization.
However, the following result proved in Part~I
states that the converse is true.
This result will play a key role in Section~\ref{gabordualsec}
for determining the properties of the
canonical dual frame of an irregular Gabor frame.

\begin{theorem} \label{selflocthm}
Let $\Fc = \set{f_i}_{i \in I}$ be a frame for $H$, and let
$a \colon I \to G$ be an associated map such that $D^+(I,a) < \infty$.
Let $\tFc$ be the canonical dual frame and
$S^{-1/2}(\Fc)$ the canonical Parseval frame.
If $(\Fc,a)$ is $\ell^1$-self-localized, then:

\begin{enumerate}
\item[(a)]
$(\Fc,a)$ is $\ell^1$-localized with respect to its
canonical dual frame $\tFc = \set{\tf_i}_{i \in I}$,

\smallskip
\item[(b)]
$(\tFc,a)$ is $\ell^1$-self-localized, and

\smallskip
\item[(c)]
($S^{-1/2}(\Fc),a)$ is $\ell^1$-self-localized.
\end{enumerate}
\end{theorem}

\smallskip
\subsection{Relative Measure} \label{relativemeasure}

\begin{definition} \label{redundancydef}
(a) Let $\Fc = \set{f_i}_{i \in I}$ and $\Ec = \set{e_j}_{j \in G}$
be frame sequences in $H$, and let 
$a \colon I \to G$ be an associated map.
Let $P_\Fc$, $P_\Ec$ denote the orthogonal projections of $H$ onto
$\clspan(\Fc)$ and $\clspan(\Ec)$, respectively.
Then given a free ultrafilter $p$ and a sequence of centers
$c = (c_N)_{N \in \N}$ in $G$, we define the
\emph{relative measure of $\Fc$ with respect to $\Ec$, $p$, and $c$} to be
$$\cM_\Ec(\Fc; p, c)
\EQ \plim_{N \in \N} \frac1{|I_N(c_N)|} \sum_{i \in I_N(c_N)}
    \ip{P_\Ec f_i}{\tf_i}.$$
The \emph{relative measure of $\Ec$ with respect to $\Fc$, $p$, and $c$} is
$$\cM_\Fc(\Ec; p, c)
\EQ \plim_{N \in \N} \frac1{|S_N(c_N)|} \sum_{j \in S_N(c_N)}
      \ip{P_\Fc \te_j}{e_j}.
$$

\smallskip
(b) If $\clspan(\Ec) \supset \clspan(\Fc)$ then $P_\Ec$ is the identity map
and $\Ec$ plays no role in determining the value of $\cM_\Ec(\Fc;p,e)$.
Therefore, in this case we define the
\emph{measure of $\Fc$ with respect to $p$ and $c$} to be
$$\cM(\Fc; p, c)
\EQ \plim_{N \in \N} \frac1{|I_N(c_N)|} \sum_{i \in I_N(c_N)}
    \ip{f_i}{\tf_i}.$$
Since $\ip{f_i}{\tf_i} = \norm{S^{-1/2}f_i}^2$, we have that
$\cM(\Fc; p, c)$ is real.
Additionally, since $S^{-1/2}(\Fc)$ is a Parseval frame, we have
$0 \le \ip{f_i}{\tf_i} \le 1$ for all $i$, and therefore
$$0 \LE \cM(\Fc;p,c) \LE 1.$$
For this case we further define the 
\emph{lower and upper measures of $\Fc$} by
\begin{align*}
\cM^-(\Fc)
& \EQ \liminf_{N \to \infty} \, \inf_{j \in G} \,
      \frac1{|I_N(j)|} \sum_{i \in I_N(j)}
      \ip{f_i}{\tf_i},
      \\[1 \jot]
\cM^+(\Fc)
& \EQ \limsup_{N \to \infty} \, \sup_{j \in G} \,
      \frac1{|I_N(j)|} \sum_{i \in I_N(j)}
      \ip{f_i}{\tf_i}.
\end{align*}
It can be seen that there exist free ultrafilters
$p^-$, $p^+$ and sequences of centers $c^-$, $c^+$ such that
$\cM^-(\Fc) = \cM(\Fc;p^-,c^-)$ and
$\cM^+(\Fc) = \cM(\Ec;p^+,c^+)$.

\medskip
(c) When $\clspan(\Fc) \supset \clspan(\Ec)$, we define the measures
$\cM(\Ec; p, c)$ and $\cM^\pm(\Ec)$ in an analogous manner.
~\qed
\end{definition}

\begin{example} \label{specialcases}
The following special cases show that the measure of a Riesz basis is~$1$.

\smallskip
\begin{enumerate}
\item[(a)] If $\clspan(\Ec) \supset \clspan(\Fc)$ and
$\Fc$ is a Riesz sequence then
$\ip{f_i}{\tf_i} = 1$ for every $i \in I$, so
$\cM(\Fc;p,c) = \cM^+(\Fc) = \cM^-(\Fc) = 1$.

\medskip
\item[(b)] If $\clspan(\Fc) \supset \clspan(\Ec)$ and
$\Ec$ is a Riesz sequence then
$\ip{\te_j}{e_j} = 1$ for every $j \in G$, so
$\cM(\Ec;p,c) = \cM^+(\Ec) = \cM^-(\Ec) = 1$.
~\qed
\end{enumerate}
\end{example}

\begin{example}[Lattice Gabor Systems] \label{latticegabor}
Consider a lattice Gabor frame, i.e., a frame of the form
$\Gc(g,\alpha\Z^d \times \beta\Z^d)$.
The canonical dual frame is a lattice Gabor frame of the form
$\Gc(\tg,\alpha\Z^d \times \beta\Z^d)$ for some $\tg \in L^2(\R^d)$.
By the Wexler--Raz relations, we have $\ip{g}{\tg} = (\alpha\beta)^d$
(we also derive this fact directly from our results in
Theorem~\ref{gabordensity}).
Since for all $k$, $n$ we have
$\ip{M_{\beta n} T_{\alpha k} g}{M_{\beta n} T_{\alpha k} \tg}
= \ip{g}{\tg}$,
it follows that for any free ultrafilter $p$ and sequence of centers
$c = (c_N)_{N \in \N}$ in $\alpha\Z^d \times \beta\Z^d$,
$$\cM(\Gc(g,\alpha\Z^d \times \beta\Z^d);p,c)
\EQ \cM^\pm(\Gc(g,\alpha\Z^d \times \beta\Z^d))
\EQ \ip{g}{\tg}
\EQ (\alpha\beta)^d.$$
Since we also have
$\BD^\pm(\alpha\Z^d \times \beta\Z^d) = (\alpha\beta)^{-d}$,
we conclude that
\begin{equation} \label{alphabeta}
\cM^\pm(\Gc(g,\alpha\Z^d \times \beta\Z^d))
\EQ \frac1{\BD^\mp(\alpha\Z^d \times \beta\Z^d)}. \quad\qed
\end{equation}
\end{example}

Equation~\eqref{alphabeta} is essentially a special case of the
Density--Relative Measure Theorem of Part~I,
which states that if $\Fc$, $\Ec$ are frame sequences such that
$(\Fc,a,\Ec)$ has both $\ell^2$-column and row decay and $D^+(I,a) < \infty$,
then
$$\cM_\Fc(\Ec;p,c) \EQ D(p,c) \cdot \cM_\Ec(\Fc;p,c)$$
for every free ultrafilter $p$ and choice of centers
$c = (c_N)_{N \in \N}$ in $G$.
We apply this general result to irregular Gabor frames
in Section~\ref{gabordensitysec}.

\smallskip
\section{Localization of Gabor Systems} \label{gaborsec}

In this section we will determinine the localization properties of
Gabor systems.

\smallskip
\subsection{Notation and Preliminary Observations for Gabor Systems}
\label{gaborprelims}

\subsubsection{Gabor systems and the reference system}
A generic Gabor system generated by a function $g \in L^2(\R^d)$ and a
sequence $\Lambda \subset \R^{2d}$
will be written in any of the following forms:
$$\Gc(g,\Lambda)
\EQ \set{M_\omega T_x g}_{(x,\omega) \in \Lambda}
\EQ \set{e^{2\pi i \omega \cdot t} g(t-x)}_{(x,\omega) \in \Lambda}
\EQ \set{g_\lambda}_{\lambda \in \Lambda}.$$
In the case that $\Gc(g,\Lambda)$ is a frame sequence we let
$$\tGc \EQ \set{\tg_\lambda}_{\lambda \in \Lambda}$$
denote the canonical dual frame sequence in $\clspan(\Gc(g,\Lambda))$,
but it is important to note that while $g_\lambda$ is a time-frequency
shift of $g$, it need not be the case that the functions $\tg_\lambda$
are time-frequency shifts of a single function.
We address the question of the structure of the dual frame in more detail
in Section~\ref{gabordualsec}.

Our reference systems will be lattice Gabor systems indexed by the group
$$G \EQ \alpha\Z^d \times \beta \Z^d,$$
where $\alpha$, $\beta > 0$ are fixed scalars.
For compactness of notation, we usually let $G$ implicitly denote the group
above, only writing out $\alpha\Z^d \times \beta\Z^d$ when
we wish to explicitly emphasize the values of $\alpha$, $\beta$.
Thus our reference systems have the form
$$\Gc(\phi, G)
\EQ \Gc(\phi, \alpha\Z^d \times \beta\Z^d)
\EQ \set{M_{\eta} T_{u} \phi}_{(\eta,u) \in G}
\EQ \set{M_{\beta n} T_{\alpha k} \phi}_{k,n \in \Z^d}.$$
The canonical dual frame of a lattice Gabor frame sequence
is another lattice Gabor frame sequence $\Gc(\tilde\phi, G)$,
generated by some dual window $\tilde\phi \in L^2(\R^d)$.
Usually the reference system makes an appearance only during the course of
a proof, and does not appear in the statement of most of the theorems.

\subsubsection{Cubes and the $a$ mapping}
For simplicity of notation, we introduce the following abbreviations.

Given $x = (x_1,\dots,x_d) \in \R^d$, we define
$\Int(x) = (\lfloor x_1 \rfloor, \dots, \lfloor x_d \rfloor)$
and $\Frac(x) = x - \Int(x)$.

Let $\alpha$, $\beta$ be fixed.
Then given a point $x = (x_1,\dots,x_d) \in \R^d$, we set
\begin{equation} \label{abbreviation}
\bx     \EQ \alpha \, \Int\bigparen{\tfrac{x}{\alpha}}, \quad
\tx     \EQ \alpha \, \Frac\bigparen{\tfrac{x}{\alpha}}, \quad
\bomega \EQ \beta \, \Int\bigparen{\tfrac{\omega}{\beta}}, \quad
\tomega \EQ \beta \, \Frac\bigparen{\tfrac{\omega}{\beta}}.
\end{equation}
Note the implicit dependence on $\alpha$ and $\beta$ in this notation.

We define the map $a \colon \Lambda \to G$ as in
Example~\ref{gaborexample} by rounding off to a near element
of $G$, i.e.,
$$a(x,\omega)
\EQ \bigparen{\alpha\Int\bigparen{\tfrac{x}{\alpha}},
              \beta\Int\bigparen{\tfrac{\omega}{\beta}}}
\EQ \bigparen{\bx, \bomega},
\qquad (x,\omega) \in \Lambda.$$

Given $z = (x,y) \in \R^{2d}$, let $Q_r(z) = Q_r(x,y)$ denote the closed
cube in $\R^{2d}$ centered at $z$ with side length $r$.
Then given $j \in G = \alpha\Z^d \times \beta\Z^d$, we have
$$S_N(j) \EQ G \cap Q_N(j)
\qquad\text{and}\qquad
I_N(j) \EQ a^{-1}(G \cap Q_N(j)).$$
Note that $I_N(j)$ is very nearly
$\Lambda \cap Q_N(j)$,
except for the effect of rounding off points via the~$a$ map.
Thus
\begin{align}
|S_N(j)|
& \EQ |G \cap Q_N(j)| \; \approx \; (\alpha\beta)^{-d} \, N^{2d},
      \label{approxnum} \\
|I_N(j)|
& \EQ |a^{-1}(G \cap Q_N(j))| \; \approx \; |\Lambda \cap Q_N(j)|. \notag
\end{align}

\subsubsection{Density and Measure} \label{beurlingsection}
Recall from \eqref{beurlingdef} the definitions of the lower
and upper Beurling densities of the index set~$\Lambda$:
$$\BD^-(\Lambda)
\EQ \liminf_{N \to \infty} \inf_{j \in \R^{2d}}
    \frac{|\Lambda \cap Q_N(j)|}{N^{2d}}
\quad\text{and}\quad
\BD^+(\Lambda)
\EQ \limsup_{N \to \infty} \sup_{j \in \R^{2d}}
    \frac{|\Lambda \cap Q_N(j)|}{N^{2d}}.$$
Note that taking the inf and sup over $j \in G$ instead of $j \in \R^{2d}$
does not affect the value of these densities.
Example~\ref{gabordensityrel} derived the relationship between
the Beurling densities $\BD^\pm(\Lambda)$ and the densities
$D^\pm(\Lambda,a)$ defined in this paper.
Specifically,
$$\BD^+(\Lambda)
\EQ (\alpha\beta)^{-d} \, D^+(\Lambda,a)
\EQ (\alpha\beta)^{-d} \,
    \limsup_{N \to \infty} \sup_{j \in \alpha\Z^d \times \beta\Z^d}
    \frac{|a^{-1}(G \cap Q_N(j))|}{|G \cap Q_N(j)|},$$ 
and similarly $\BD^-(\Lambda) = (\alpha\beta)^{-d} \, D^-(\Lambda,a)$.
In light of this equation, we define the Beurling density of $\Lambda$
with respect to a free ultrafilter~$p$ and a sequence of centers
$c = (c_N)_{N \in \N}$ in $\R^{2d}$ to be
$$\BD(\Lambda;p,c)
\EQ (\alpha\beta)^{-d} \, D(\Lambda,a;p,c)
\EQ (\alpha\beta)^{-d} \, \plim_{N \in \N}
    \frac{|a^{-1}(G \cap Q_N(c_N))|} {|G \cap Q_N(c_N)|}.$$
Our results for Gabor systems will all be stated in terms of these
Beurling densities.

The measure of a Gabor frame sequence $\Gc(g,\Lambda)$
with respect to a free ultrafilter $p$ and a sequence
of centers $c = (c_N)_{N \in \N}$ in $\R^{2d}$ is
$$\cM(\Gc(g,\Lambda); p, c)
\EQ \plim_{N \in \N} \frac1{|a^{-1}(G \cap Q_N(c_N))|}
    \sum_{a^{-1}(G \cap Q_N(c_N))} \ip{g_\lambda}{\tg_\lambda}.$$

In particular, note that if $\Gc(g,\Lambda)$ is a Riesz sequence then
$\ip{g_\lambda}{\tg_\lambda} = 1$ for all $\lambda$, so the measure is
$\cM^\pm(\Gc(g,\Lambda)) = 1$.  
Also recall that the measure of a lattice Gabor system was computed in
Example~\ref{latticegabor}.

By making the approximations in \eqref{approxnum} precise, the next
lemma reformulates the density and measure in such a way that it becomes
clear that they do not depend on the choice of $\alpha$, $\beta$.
That is, the density of $I$ and the measure of $\Gc(g,\Lambda)$
are independent of the choice of reference group.
Analogous reformulations of the upper and lower density and measures also
hold under the same hypotheses.

\begin{lemma} \label{reformulation}
Let $\Lambda \subset \R^{2d}$ be given.

\begin{enumerate}
\item[(a)]
If $\BD^+(\Lambda) < \infty$, then for any ultrafilter $p$ and any
sequence of centers $c = (c_N)_{N \in \N}$ in~$\R^{2d}$,
$$\BD(\Lambda;p,c)
\EQ \plim_{N \in \N} \frac{|\Lambda \cap Q_N(c_N)|} {N^{2d}}.$$

\item[(b)]
Let $g \in L^2(\R^d)$ be given.
If $0 < \BD^-(\Lambda) \le \BD^+(\Lambda) < \infty$,
then for any ultrafilter $p$ and any
sequence of centers $c = (c_N)_{N \in \N}$ in~$\R^{2d}$,
$$\cM(\Gc(g,\Lambda);p,c)
\EQ \plim_{N \in \N}
    \frac1{|\Lambda \cap Q_N(c_N)|} \sum_{\lambda \in \Lambda \cap Q_N(c_N)} \,
    \ip{g_\lambda}{\tg_\lambda}.$$
\end{enumerate}
\end{lemma}

\begin{proof}
(a) The map $a$ is a bounded perturbation of the identity map.
In particular, if $\delta = \max\set{\alpha,\beta}$, then we have
\begin{equation} \label{inclusions}
\Lambda \cap Q_{N-\delta}(c_N)
\SUBSET a^{-1}(G \cap Q_N(c_N))
\SUBSET \Lambda \cap Q_{N+\delta}(c_N).
\end{equation}
Since the upper density is finite, there is a constant $C$ such that
$|\Lambda \cap Q_{N+\delta}(c_N) \setminus Q_N(c_n)| \le C \, N^{2d-1}$.
Using the second inclusion in \eqref{inclusions}, we therefore have
\begin{align*}
\BD(\Lambda;p,c)
& \EQ (\alpha\beta)^{-d} \, \plim_{N \in \N}
      \frac{|a^{-1}(G \cap Q_N(c_N))|} {|G \cap Q_N(c_N)|} \\[1 \jot]
& \LE (\alpha\beta)^{-d} \, \plim_{N \in \N}
      \frac{|\Lambda \cap Q_{N+\delta}(c_N)|} {|G \cap Q_N(c_N)|}
      \allowdisplaybreaks \\[1 \jot]
& \EQ (\alpha\beta)^{-d} \, \plim_{N \in \N}
      \biggparen{\frac{|\Lambda \cap Q_N(c_N)|} {|G \cap Q_N(c_N)|} \plus
      \frac{|\Lambda \cap Q_{N+\delta}(c_N) \setminus Q_N(c_N)|}
           {|G \cap Q_N(c_N)|}} \\[1 \jot]
& \EQ \plim_{N \in \N}
      \frac{|\Lambda \cap Q_N(c_N)|} {N^{2d}},
\end{align*} 
and a similar computation making use of the first inclusion
in \eqref{inclusions} yields the opposite inequality.

\medskip
(b) The fact that $0 < \BD^-(\Lambda)$ and $\BD^+(\Lambda) < \infty$
implies that there exist $C_1$, $C_2 > 0$ such that for all $N$ large
enough and all $j \in G$ we have
$C_1 N^{2d} \le |\Lambda \cap Q_N(j)| \le C_2 N^{2d}$.
Combining this with \eqref{inclusions} then implies that
$$\plim_{N \in \N} \frac{|\Lambda \cap Q_N(c_N)|}{|a^{-1}(G \cap Q_N(c_N))|}
\EQ 1.$$
Since $0 \le \ip{g_\lambda}{\tg_\lambda} \le 1$, we therefore have
\begin{align*}
\cM(\Gc(g,\Lambda);p,c)
& \EQ \plim_{N \in \N} \frac1{|a^{-1}(G \cap Q_N(c_N))|}
      \sum_{\lambda \in a^{-1}(G \cap Q_N(c_N))} \,
      \ip{g_\lambda}{\tg_\lambda} \\[1 \jot]
& \LE \plim_{N \in \N} \frac1{|\Lambda \cap Q_N(c_N)|}
      \sum_{\lambda \in \Lambda \cap Q_{N+\delta}(c_N)} \,
      \ip{g_\lambda}{\tg_\lambda}
      \allowdisplaybreaks \\[1 \jot]
& \LE \plim_{N \in \N} \frac1{|\Lambda \cap Q_N(c_N)|}
      \sum_{\lambda \in \Lambda \cap Q_N(c_N)} \,
      \ip{g_\lambda}{\tg_\lambda} \plus \\
& \qquad\quad \plim_{N \in \N}
      \frac{|\Lambda \cap Q_{N+\delta}(c_N) \setminus Q_N(c_N)|}
           {|\Lambda \cap Q_N(c_N)|},
\end{align*}
and the final term in this computation is zero.
Combining this with a similar computation for the opposite inequality then
yields the result.  
\end{proof}

\smallskip
\subsection{Localization with respect to the Box Function} \label{boxlocalsec}

We now show that any Gabor system $\Gc(g,\Lambda)$
has $\ell^2$-row decay with respect to
the Gabor orthonormal basis $\Gc(\CHI,\Z^{2d})$ generated by the
``box function'' $\CHI = \CHI_{[-\frac12,\frac12)^d}$,
and that we recover as a consequence the fact first proved in \cite{CDH99}
that any Gabor system that forms a Bessel sequence
must have finite density.

The Gabor system generated by the box function is convenient both because
it is an orthonormal basis and because
the index set is $G = \Z^{2d}$ (so, in particular, $\alpha = \beta = 1$).
However, in general this is not a useful basis in applications because
the generating function~$\CHI$ is poorly concentrated in the
time-frequency plane (in fact, by the Balian--Low Theorem, no Gabor Riesz
basis of the form $\Gc(\phi,\alpha\Z^d \times \beta\Z^d)$
can have a generator $\phi$ that is simultaneously well-concentrated in both
time and frequency).
In the next section we show in more detail how time-frequency concentration
is related to localization.

\begin{proposition} \label{boxlocal}
If $g \in L^2(\R^d)$ and $\Lambda \subset \R^{2d}$, then
$(\Gc(g,\Lambda), \, a, \, \Gc(\CHI,\Z^{2d}))$
has $\ell^2$-row decay.
\end{proposition}
\begin{proof}
Choose $\eps > 0$, and let $R>0$ be large enough that
$\int_{\R^d \setminus [-R,R]^d} |g(t)|^2 \, dt < \eps$.
Fix an even integer $N_\eps > R+3$.

Consider now any $(x,\omega) \in \Lambda$.
Since $\alpha = \beta = 1$, we have
$\tx = \Frac(x) \in [0,1)^d$ and
$\tomega = \Frac(\omega) \in [0,1)^d$.
Note that for each $k \in \Z^d$,
$\set{M_n T_k \CHI}_{n \in \Z^d}$ is an orthonormal basis for
the subspace of $L^2(\R^d)$ consisting of functions supported in the
unit cube $B_k$ in $\R^d$ centered at~$k$.
Therefore,
\begin{align*}
\sum_{(k,n) \in \Z^{2d} \setminus S_{N_\eps}(a(x,\omega))}
|\ip{M_\omega T_x g} {M_n T_k \CHI}|^2
& \EQ \sum_{(k,n) \in \Z^{2d} \setminus S_{N_\eps}(\bx, \bomega)}
      |\ip{M_\tomega T_\tx g} {M_{n-\bomega} T_{k-\bx} \CHI}|^2
      \allowdisplaybreaks \\[1 \jot]
& \EQ \sum_{(k,n) \in \Z^{2d} \setminus S_{N_\eps}(0,0)}
      |\ip{M_\tomega T_\tx g} {M_n T_k \CHI}|^2
      \allowdisplaybreaks \\[1 \jot]
& \LE \sum_{k \in \Z^d \setminus [-\frac{N_\eps}2,\frac{N_\eps}2)^d} \,
      \sum_{n \in \Z^d} |\ip{M_\tomega T_\tx g} {M_n T_k \CHI}|^2
      \allowdisplaybreaks \\[1 \jot]
& \EQ \sum_{k \in \Z^d \setminus [-\frac{N_\eps}2,\frac{N_\eps}2)^d} \,
      \int_{B_k} |M_\tomega T_\tx g(t)|^2 \, dt
      \allowdisplaybreaks \\[1 \jot]
& \LE \int_{\R^d \setminus [\frac{-N_\eps+1}2,\frac{N_\eps-1}2]^d}
      |g(t-\tx)|^2 \, dt
      \\[1 \jot]
& \LE \int_{\R^d \setminus [-\frac{R}2,\frac{R}2]^d} |g(t)|^2 \, dt
\LT \eps. \qedhere
\end{align*}
\end{proof}

Note that the preceeding result does not require that $\Lambda$ have
finite density.
However, we next observe that it follows as a consequence of Part~I
results that the density must be finite if the Gabor system
is a Bessel sequence.

\begin{corollary} \label{gaborfinite}
If $g \in L^2(\R^d)$ and $\Lambda \subset \R^{2d}$
are such that $\Gc(g,\Lambda)$
is a Bessel sequence, then $\BD^+(\Lambda) < \infty$.
\end{corollary}
\begin{proof}
Since $(\Gc(g,\Lambda), \, a, \, \Gc(\CHI,\Z^{2d}))$
has $\ell^2$-row decay and $\Gc(g,\Lambda)$ is norm-bounded below,
all the hypotheses of \cite[Thm.~3.3]{BCHL05a} are fulfilled,
and consequently the upper density must be finite.
\end{proof}

\smallskip
\subsection{Localization with respect to $M^1$ Functions} \label{M1localsec}

For most applications in time-fre\-quency analysis, the generator of a Gabor
system must possess some amount of joint concentration in both time
and frequency.
Concentration is quantified by the norms of the modulation spaces, which
are the Banach function spaces naturally associated to time-frequency
analysis.
The modulation spaces were invented and extensively investigated by
Feichtinger, with some of the main references being
\cite{Fei81}, \cite{Fei89}, \cite{FG89a}, \cite{FG89b}, \cite{FG97},
\cite{Fei03}.
For a detailed development of the theory of modulation spaces and their
weighted counterparts, we refer to the original literature mentioned above and
to \cite[Ch.~11--13]{Gro01}.

For our purposes, the following special case of unweighted modulation
spaces will be sufficient.

\begin{definition}\label{gaussian} \
\begin{enumerate}
\item[(a)]
The \emph{Short-Time Fourier Transform} (STFT)
of a tempered distribution $g \in \Sc'(\R^d)$ with
respect to a window function $\phi \in \Sc(\R^d)$ is
$$V_\phi g(x,\omega) \EQ \ip{g}{M_\omega T_x \phi},
\qquad (x,\omega) \in \R^{2d}.$$

\medskip
\item[(b)]
Let $\gamma(x) = 2^{d/4} e^{-\pi x \cdot x}$ be the Gaussian
function, which has been normalized so that $\norm{\gamma}_2 = 1$.
Then for $1 \le p \le \infty$, the modulation space $M^p(\R^d)$ consists of
all tempered distributions $f \in \Sc'(\R^d)$ such that
\begin{equation} \label{mpnorm}
\norm{f}_{M^p}
\EQ \norm{V_\gamma f}_{L^p}
\EQ \biggparen{\iint_{\R^{2d}} |\ip{f}{M_\omega T_x \gamma}|^p \,
    dx \, d\omega}^{1/p}
\LT \infty,
\end{equation}
with the usual adjustment if $p=\infty$.
~\qed
\end{enumerate}
\end{definition}

$M^p$ is a Banach space for each $1 \le p \le \infty$,
and any nonzero function $g \in M^1$
can be substituted for~$\gamma$ in \eqref{mpnorm} to define an
equivalent norm for $M^p$.
We have $M^2 = L^2$, and
$\Sc \subsetneq M^p \subsetneq M^q \subsetneq \Sc'$
for $1 \le p < q \le \infty$,
where $\Sc$ is the Schwartz class.
The box function $\CHI$ lies in $M^p$ for $p>1$, but is not in $M^1$.

\begin{remark}
There is a complete characterization of the frame properties of lattice
Gabor systems generated by the Gaussian:

\smallskip
\begin{enumerate}
\item[(a)] $\Gc(\gamma,\alpha\Z^d \times \beta\Z^d)$
is a frame for $L^2(\R^d)$ if $0<\alpha\beta<1$,

\item[(b)] $\Gc(\gamma,\alpha\Z^d \times \beta\Z^d)$
is a Riesz sequence in $L^2(\R^d)$ if $\alpha\beta>1$
(but not a Riesz basis),

\item[(c)] $\Gc(\gamma,\alpha\Z^d \times \beta\Z^d)$
is complete but not a frame for $L^2(\R^d)$ if $\alpha\beta=1$.
\end{enumerate}

\smallskip\noindent
Part~(a) was proved in \cite{Lyu92}, \cite{Sei92}, \cite{SW92};
see also the simple proof given in \cite{Jan94}.
Part~(b) is a consequence of part~(a) and the Wexler--Raz relations, and
part~(c) is easy to show using Zak transform techniques, see \cite{Gro01}.
Such a complete characterization of frame properties is known for only a
few specific functions \cite{Jan94}, \cite{Jan03}, \cite{JS02},
cf.\ also \cite{CK02}.
In particular such a characterization is not
available for general $M^1$ functions.
On the other hand, given any particular choice of $\alpha$, $\beta$ with
$0<\alpha\beta<1$, it is easy to construct a function
$g \in C_c^\infty(\R^d)$ such that $\Gc(g,\alpha\Z^d \times \beta\Z^d)$
is a Parseval frame for $L^2(\R^d)$, cf.\ \cite{DGM86}.
It is known that if $\Gc(\phi, \alpha\Z^d \times \beta\Z^d)$
is a frame for $L^2(\R^d)$ generated by a function $\phi \in M^1$,
then the canonical dual frame
$\Gc(\tilde\phi, \alpha\Z^d \times \beta\Z^d)$
has a generator $\tilde\phi$ that also lies in $M^1$ \cite{GL04}.
In Section~\ref{gabordualsec} we will extend that result to the more
general setting of irregular Gabor frames.
~\qed
\end{remark}

In addition to the modulation spaces, we will also need a special
case of the Wiener amalgam spaces on $\R^{2d}$.
Feichtinger has developed a general notion of amalgam spaces, e.g.,
\cite{Fei80}, \cite{FG85}, \cite{Fei87}, \cite{Fei90}, \cite{Fei92}.
For an introduction, with extensive references to the original literature,
to the particular Wiener amalgams appearing in the
following definition, we refer to \cite{Hei03}.

\begin{definition} \label{amalgamdef}
Given $1 \le p \le \infty$, the Wiener amalgam $W(\Cc,\ell^p)$ consists
of all continuous functions $F$ on $\R^{2d}$ for which
$$\norm{F}_{W(\Cc,\ell^p)}
\EQ \biggparen{\sum_{(k,n) \in \Z^{2d}} \,
    \sup_{(u,\eta) \in \Q_{\alpha,\beta}(\alpha k, \beta n)} \,
    |F(u,\eta)|^p}^{1/p}
\LT \infty,$$
with the usual adjustment if $p=\infty$, and where
$\Q_{\alpha,\beta}(x,y) = [0,\alpha)^d \times [0,\beta)^d + (x,y)$.
~\qed
\end{definition}

$W(\Cc,\ell^p)$ is a Banach space, and its definition is independent of
the values of $\alpha$ and $\beta$ in the sense that each choice of
$\alpha$, $\beta$ yields an equivalent norm for $W(\Cc,\ell^p)$.

We will require the following lemma on the basic properties of the STFT.
Part~(a) is proved in \cite[Thm.~12.2.1]{Gro01},
and part~(b) in \cite[Lem12.1.1]{Gro01}.

\begin{lemma} \label{pointest} \

\begin{enumerate}
\item[(a)]
Let $1 \le p \le \infty$ be given.
If $g \in M^p$ and $\phi \in M^1$, then $V_\phi g \in W(\Cc,\ell^p)$, and
$$\norm{V_\phi g}_{W(\Cc,\ell^p)}
\LE C \, \norm{g}_{M^p} \, \norm{\phi}_{M^1},$$
where $C$ is a constant independent of $g$ and $\phi$.

\medskip
\item[(b)]
Let $f$, $g\in L^2(\R^d)$ be given.
If $V_g f \in L^1(\R^{2d})$, then $f$, $g \in M^1(\R^d)$.
\end{enumerate}
\end{lemma}

Next we will show that if the generator of our reference system is an
$M^1$ function $\phi$, then
$(\Gc(g,\Lambda), \, a, \, \Gc(\phi, G))$
is $\ell^p$-localized whenever $g \in M^p$.
We also show that the converse is true if we assume that $\phi$
generates a frame.
Thus the degree of localization is tied to the time-frequency concentration
of the generator $g$.
This is a stronger statement than previously known results, which only
demonstrated that Gabor frames satisfy the weak HAP with respect to
reference systems generated by the Gaussian.
We will also show that if $g \in M^1$ then $(\Gc(g,\Lambda),a)$ is
$\ell^1$-self-localized (this statement does not require a reference system).
However, the converse of this is false.
For example, if we set $g = \CHI$
and $\Lambda = \Z^{2d}$, then $\Gc(g,\Lambda)$ is an orthonormal basis for
$L^2(\R^d)$ and hence is $\ell^1$-self-localized, but $g \notin M^1$.

\begin{theorem} \label{M1local}
Let $g \in L^2(\R^d)$ and $\Lambda \subset \R^{2d}$ be given.
Let $\phi \in L^2(\R^d)$ and $\alpha$, $\beta > 0$ be given, and
fix $1 \le p \le 2$.
Then the following statements hold.

\smallskip
\begin{enumerate}
\item[(a)] If $g \in M^p$ and $\phi \in M^1$ then
$(\Gc(g,\Lambda), \, a, \, \Gc(\phi, G))$
is $\ell^p$-localized.

\medskip
\item[(b)]
Suppose $\phi \in M^1$ and $\alpha$, $\beta > 0$ are such that
$\Gc(\phi,G)$ is a frame for $L^2(\R^d)$.
If $(\Gc(g,\Lambda), \, a, \, \Gc(\phi, G))$
is $\ell^p$-localized, then $g \in M^p$.

\medskip
\item[(c)] If $g \in M^1$ and $\phi \in M^p$ then
$(\Gc(g,\Lambda), \, a, \, \Gc(\phi, G))$
is $\ell^p$-localized.

\medskip
\item[(d)]
If $g \in M^1$ then
$(\Gc(g,\Lambda),a)$ is $\ell^1$-self-localized.
\end{enumerate}
\end{theorem}
\begin{proof}
(a) Set
$$r_{(\alpha k,\beta n)}
\EQ \sup_{(u,\eta) \in \Q_{\alpha,\beta}(\alpha k, \beta n)} \,
    |V_\phi g(-u,-\eta)|.$$
By Lemma~\ref{pointest} we have $V_\phi g \in W(\Cc,\ell^p)$, so
$r = (r_{(\alpha k,\beta n)})_{(\alpha k,\beta n) \in G} \in \ell^p(G)$.
Let $(x,\omega) \in \Lambda$ and $(u,v) \in G$ be given.
Then, recalling the notations~$\tx$, $\tomega$ introduced
in \eqref{abbreviation},
since $\tx \in [0,\alpha)^d$ and $\tomega \in [0,\beta)^d$, we have
\begin{align*}
|\ip{M_\omega T_x g}{M_{v} T_{u} \phi}|
& \EQ |\ip{g} {M_{v - \bomega - \tomega} T_{u - \bx - \tx} \phi}| \\[1 \jot]
& \EQ |V_\phi g(u - \bx - \tx, v - \bomega - \tomega)|
\LE r_{a(x,\omega) - (u,v)},
\end{align*}
so $\ell^p$-localization holds.

\medskip
(b) By \cite{GL04} or Theorem~\ref{dualgaborlocalization}, the dual frame
of $\Gc(\phi,G)$ is a lattice Gabor frame $\Gc(\tilde\phi,G)$
with a dual window $\tilde\phi \in M^1$.
Fix any $(x,\omega) \in \Lambda$.
Expanding $M_\omega T_x g$ with respect to this frame, we have
\begin{equation} \label{gexpand}
M_\omega T_x g
\EQ \sum_{k,n \in \Z^d} \ip{M_\omega T_x g}{M_{\beta n} T_{\alpha k} \phi} \,
    M_{\beta n} T_{\alpha k} \tilde\phi,
\end{equation}
with convergence in $L^2(\R^d)$.
Now, by definition of $\ell^p$-localization, there exists $r \in \ell^p(G)$
such that
$$|\ip{M_\omega T_x g}{M_{\beta n} T_{\alpha k} \phi}|
\LE r_{a(x,\omega) - (\alpha k, \beta n)},
\qquad (\alpha k, \beta n) \in G.$$
Consequently,
$\set{\ip{M_\omega T_x g}{M_{\beta n} T_{\alpha k} \phi}}_{k,n \in \Z^d}
\in \ell^p$,
and so by \cite[Thm.~12.2.4]{Gro01} the series on the right-hand side
of \eqref{gexpand} converges in $M^p$-norm.
Since it also converges in $L^2$-norm, the series must converge in $M^p$-norm
to $M_\omega T_x g$, so $M_\omega T_x g \in M^p$.
Since $M^p$ is closed under time-frequency shifts, we conclude $g \in M^p$.

\medskip
(c) Since $|V_\phi g(x,\omega)| = |V_g \phi(-x,-\omega)|$,
this follows from part~(a).

\medskip
(d) This can be shown directly, similarly to part~(a), or by applying
\cite[Lem.~2.15]{BCHL05a} to part~(a), using as a reference system any
lattice Gabor frame $\Gc(\phi, G)$ whose generator $\phi$ lies in $M^1$.
\end{proof}

\smallskip
\section{New Implications for Gabor Frames} \label{newresults}

\smallskip
\subsection{Density and Overcompleteness for Gabor Systems}
\label{gabordensitysec}
Parts~(b) and~(c) of the following theorem are new results for Gabor frames; the equalities given are much stronger than the single inequality obtained in \cite{BCHL03b}.
Parts~(a) and~(d)  recover the known density facts for
irregular Gabor frames, and part~(e) is the special case of lattice systems.
These results are stated in terms of Beurling density, which was discussed
in Section~\ref{beurlingsection}.

\begin{theorem} \label{gabordensity}
Let $g \in L^2(\R^d)$ and $\Lambda \subset \R^{2d}$ be such that
$\Gc(g,\Lambda)$ is a Gabor frame for $L^2(\R^d)$.
Then the following statements hold.

\smallskip
\begin{enumerate}
\item[(a)]
$1 \le \BD^-(\Lambda) \le \BD^+(\Lambda) < \infty$.

\medskip
\item[(b)]
For any free ultrafilter $p$ and any
sequence of centers $c = (c_N)_{N \in \N}$ in $\R^d$, we have
\begin{equation} \label{reciprocal}
\cM(\Gc(g,\Lambda); p, c) \EQ \frac{1}{\BD(\Lambda; p, c)},
\end{equation}
and consequently
\begin{equation} \label{averaging}
\plim_{N \in \N} \frac1{N^{2d}}
\sum_{\lambda \in \Lambda \cap Q_N(c_N)} \ip{g_\lambda}{\tg_\lambda}
\EQ 1.
\end{equation}

\medskip
\item[(c)]
$\cM^-(\Gc(g,\Lambda)) \EQ \dfrac1{\BD^+(\Lambda)}$
and
$\cM^+(\Gc(g,\Lambda)) \EQ \dfrac1{\BD^-(\Lambda)}$.

\medskip
\item[(d)]
If $\Gc(g,\Lambda)$ is a Riesz basis, then
$\BD^-(\Lambda) = \BD^+(\Lambda) = 1$.

\bigskip
\item[(e)]
If $\Lambda = \alpha\Z^d \times \beta\Z^d$ then
$0 < \alpha\beta \le 1$ and
$\ip{g}{\tg} = (\alpha\beta)^d$.
\end{enumerate}
\end{theorem}

\begin{proof}
(a) In this part we use a reference system
$\Gc(\gamma, G)$ generated by the Gaussian $\gamma$.
If we take any $\alpha$, $\beta > 0$ so that $\alpha\beta > 1$,
then $\Gc(\gamma, G)$ is a Riesz sequence in $L^2(\R^d)$, and
$(\Gc(g,\Lambda), \, a, \, \Gc(\gamma, G))$
is $\ell^2$-localized by Theorem~\ref{M1local}.
Therefore, we have by \cite[Thm.~3.2]{BCHL05a} that
$(\alpha\beta)^d \, \BD^-(\Lambda) = D^-(\Lambda,a) \ge 1$.
Since this is true for any $\alpha\beta > 1$,
we conclude that $\BD^-(\Lambda) \ge 1$.
The fact that $\BD^+(\Lambda) < \infty$ follows from
Corollary~\ref{gaborfinite}.

\medskip
(b) In this part our reference system will be generated by the function
$$\phi(x)
\EQ \prod_{k=1}^d \frac{e^{2\pi i x_k} + 1}2 \,
    \CHI_{[-\frac12,\frac12]}(x_k),
\qquad x = (x_1,\ldots,x_d) \in \R^d.$$
Since $\phi$ is compactly supported,
$\phi \in L^1(\R^d)$, and
$\hat\phi(\omega)
= \prod_{k=1}^d \frac{\sin \pi\omega_k}{2\pi(\omega_k - \omega_k^2)}
\in L^1(\R^d)$,
it follows from \cite[Thm.~3.2.17]{FZ98} that $\phi \in M^1$.
If we set $\alpha = 1/2$ and $\beta = 1$, then since
$\sum |\phi(x - \frac12 k)|^2 = 1$, we have by
\cite[Thm.~6.4.1]{Gro01} that
$\Gc(\phi, G)$ is a Parseval frame for $L^2(\R^d)$.
By direct computation,
the measure of this frame is
$\cM(\Gc(\phi, G); p, c)
= \norm{\phi}_2^2 = 2^{-d}$.

Since $(\Gc(g,\Lambda), \, a, \, \Gc(\phi, G))$
is $\ell^2$-localized, we have by \cite[Thm.~3.5(a)]{BCHL05a} that
$$\cM(\Gc(g,\Lambda); p, c)
\EQ \frac{\cM(\Gc(\phi, G);p,c)} {D(p,c)}
\EQ \frac{2^{-d}}{2^{-d} \, \BD(\Lambda; p, c)}
\EQ \frac1{\BD(\Lambda; p, c)}.$$

Finally, the reformulation in \eqref{averaging} follows by applying
Lemma~\ref{reformulation} to \eqref{reciprocal}.

\medskip
(c) Using the same reference system as in part~(b), we have
by \cite[Thm.~3.5(a)]{BCHL05a} that
$$2^{-d}
\EQ \cM^-(\Gc(\phi, G))
\LE D^+(\Lambda,a) \cdot \cM^-(\Gc(g,\Lambda))
\LE \cM^+(\Gc(\phi, G))
\EQ 2^{-d},$$
so
$\cM^-(\Gc(g,\Lambda))
= \frac{2^{-d}}{D^+(\Lambda,a)}
= \frac1{\BD^+(\Lambda)}$,
and similarly $\cM^+(\Gc(g,\Lambda)) = \frac1{\BD^-(\Lambda)}$.

\medskip
(d) If $\Gc(g,\Lambda)$ is a Riesz basis then
$\cM^\pm(\Gc(g,\Lambda)) = 1$, so the result follows from part~(c).

\medskip
(e) In the lattice case $\Lambda = \alpha\Z^d \times \beta\Z^d$ we have
$\BD^\pm(\alpha\Z^d \times \beta\Z^d) = (\alpha\beta)^{-d}$
and $\cM^\pm(\Gc(g,\Lambda)) = \ip{g}{\tg}$,
so the result follows from parts~(a) and~(c).
\end{proof}

Next we prove results on the relationship between the density, frame
bounds, and norm of the generator of a Gabor frame.
The special case of lattice systems was first proved by Daubechies
\cite[Eq.~2.2.9]{Dau90}.

\begin{theorem} \label{gaborbounds}
Let $g \in L^2(\R^d)$ and $\Lambda \subset \R^{2d}$ be such that
$\Gc(g,\Lambda)$ is a Gabor frame for $L^2(\R^d)$, with frame bounds $A$, $B$.
Then the following statements hold.

\smallskip
\begin{enumerate}
\item[(a)]
$A \LE \BD^-(\Lambda) \, \norm{g}_2^2
\LE \BD^+(\Lambda) \, \norm{g}_2^2
\LE B$.

\medskip
\item[(b)]
If $\Gc(g,\Lambda)$ is a tight frame, then $\Lambda$ has uniform
Beurling density, that is,
$\BD^-(\Lambda) = \BD^+(\Lambda)$,
and furthermore $A = \BD^\pm(\Lambda) \, \norm{g}_2^2$.

\medskip
\item[(c)]
If $\Lambda = \alpha\Z^d \times \beta\Z^d$, then
$A \LE \dfrac{\norm{g}_2^2}{(\alpha\beta)^d} \LE B$.
\end{enumerate}
\end{theorem}

\begin{proof}
For a reference system,
fix any $\phi \in M^1$ and any $\alpha$, $\beta > 0$
such that $\Gc(\phi,G)$ is a Parseval frame
for $L^2(\R^d)$.
By Theorem~\ref{gabordensity}(e)
we have $\norm{\phi}_2^2 = (\alpha\beta)^d$.
Further,
$(\Gc(g,\Lambda), \, a, \, \Gc(\phi, G))$
is $\ell^2$-localized by Theorem~\ref{M1local},
so \cite[Thm.~3.6(c)]{BCHL05a} implies that
\begin{align*}
\frac{A \, (\alpha\beta)^d}{\norm{g}_2^2}
\EQ \frac{A \, \norm{\phi}_2^2}{\norm{g}_2^2}
& \LE D^-(\Lambda,a) \\
& \EQ (\alpha\beta)^d \, \BD^-(\Lambda) \allowdisplaybreaks \\[1 \jot]
& \LE (\alpha\beta)^d \, \BD^+(\Lambda) \\
& \EQ D^+(\Lambda,a)
\EQ \frac{B \, \norm{\phi}_2^2}{\norm{g}_2^2}
\LE \frac{B \, (\alpha\beta)^d}{\norm{g}_2^2},
\end{align*}
so part~(a) follows.
Parts~(b) and~(c) are immediate consequences of part~(a).
\end{proof}

\smallskip
\subsection{Excess of Gabor Frames} \label{gaborremovalsec}

In this section we consider the excess of Gabor frames, and
show that subsets with positive density may be removed from an overcomplete
Gabor frame yet still leave a frame.

\begin{theorem} \label{gaborremoval}
Let $g \in L^2(\R^d)$ and $\Lambda \subset \R^{2d}$ be such that
$\Gc(g,\Lambda)$ is a Gabor frame for $L^2(\R^d)$.
Then the following statements hold.

\smallskip
\begin{enumerate}
\item[(a)]
If $\BD^+(\Lambda) > 1$ then $\Gc(g,\Lambda)$ has infinite excess, and
there exists an infinite subset $J \subset \Lambda$ such that
$\Gc(g,\Lambda \setminus J)$ is a frame for $L^2(\R^d)$.

\medskip
\item[(b)]
If $g \in M^1$ and $\BD^-(\Lambda) > 1$, then there exists
$J \subset \Lambda$ with $\BD^+(J) = \BD^-(J) > 0$ such that
$\Gc(g,\Lambda \setminus J)$ is a frame for $L^2(\R^d)$.
\end{enumerate}
\end{theorem}

\begin{proof}
(a) By Theorem~\ref{gabordensity}(c) we have $\cM^-(\Gc(g,\Lambda)) < 1$.
The result therefore follows from \cite[Prop.~2.21]{BCHL05a}.

\medskip
(b) By Theorem~\ref{gabordensity}(c) we have $\cM^+(\Gc(g,\Lambda)) < 1$.
Since $g \in M^1$, we have that $(\Gc(g,\Lambda),a)$
is $\ell^1$-self-localized by Theorem~\ref{M1local}.
Hence the result follows from \cite[Thm.~3.8]{BCHL05a}.
\end{proof}

\smallskip
\subsection{Localization and Structure of the Canonical Dual Frame}
\label{gabordualsec}

In this section we study the structure of the canonical dual frame of an
irregular Gabor frame.

First we introduce the notion of Gabor molecules.
The term ``molecule'' arises from the convention that the generator $g$ of a
Gabor system $\Gc(g,\Lambda)$ is often referred to as an ``atom.''

\begin{definition}
Let $\Lambda \subset \R^{2d}$ and
$f_\lambda \in L^2(\R^d)$ for $\lambda \in \Lambda$ be given.
Then $\Fc = \set{f_\lambda}_{\lambda \in \Lambda}$
is a \emph{set of Gabor molecules} if
there exists an envelope function $\Gamma \in W(\Cc,\ell^2)$ such that
$$\forall\, \lambda \in \Lambda, \quad
\forall\, z \in \R^{2d}, \quad
|V_\gamma f_\lambda(z)| \LE \Gamma(z - \lambda). \qeddeff$$
\end{definition}

Thus, if $\Gamma$ is concentrated around the origin in $\R^{2d}$, then
the STFT of $f_\lambda$ is concentrated around the point~$\lambda$.
Every Gabor system $\Gc(g,\Lambda)$ is a set of Gabor molecules, as
$|V_\gamma g_\lambda(z)| = |V_\gamma g(z-\lambda)|$ for every $z$, $\lambda$.
Sometimes the following equivalent definition is more convenient:
$\Fc = \set{M_\omega T_x f_{x\omega}}_{(x,\omega) \in \Lambda}$
is a set of Gabor molecules if
there exists $\Gamma \in W(\Cc,\ell^2)$ such that
$|V_\gamma f_{x\omega}(z)| \LE \Gamma(z)$
for all $(x,\omega) \in \Lambda$ and
$z \in \R^{2d}$.

The following lemma shows that the definition of Gabor molecules is unchanged
if we replace the Gaussian window by any window function $\phi \in M^1$.

\begin{lemma}
Suppose $\Fc = \set{f_\lambda}_{\lambda \in \Lambda}$
is a set of Gabor molecules with envelope $\Gamma \in W(\Cc,\ell^2)$.
If $\phi \in M^1$, then
$\Gamma_\phi = \Gamma * V_\phi \gamma \in W(\Cc,\ell^2)$, and
$|V_\phi f_\lambda(z)| \le \Gamma_\phi(z-\lambda)$
for all $\lambda \in \Lambda$ and $z \in \R^{2d}$.
\end{lemma}
\begin{proof}
Since $\gamma$, $\phi \in M^1$ we have $V_\gamma \gamma \in L^1(\R^{2d})$.
Therefore
$\Gamma_\phi \in W(\Cc,\ell^2) * L^1 \subset W(\Cc,\ell^2)$
by \cite[Thm.~11.1.5]{Gro01}.
Further,
\begin{align*}
|V_\phi f_\lambda(z)|
\LE (|V_\gamma f_\lambda| * |V_\phi \gamma|)(z)
& \EQ \int_{\R^{2d}} |V_\phi f_\lambda(z-w)| \, |V_\phi \gamma(w)| \, dw
      \allowdisplaybreaks \\
& \LE \int_{\R^{2d}} \Gamma(z-w-\lambda) \, |V_\phi \gamma(w)| \, dw
      \\[1 \jot]
& \EQ (\Gamma * |V_\phi \gamma|)(z-\lambda)
\EQ \Gamma_\phi(z-\lambda),
\end{align*}
the first inequality following from \cite[Lem.~11.3.3]{Gro01}.
\end{proof}

Gr\"{o}chenig and Leinert \cite{GL04} proved that if $\Lambda$ is a lattice
then the canonical dual frame of a lattice Gabor frame generated by a function
$g \in M^1$ is generated by a dual window that also lies in $M^1$
(they also obtained weighted versions of this result).
Their proof relied on deep results about symmetric Banach algebras.
Here we extend this result to the more general setting of
irregular Gabor frame sequences,
and furthermore determine the structure of the dual frame (which in the
lattice setting is simply another lattice Gabor frame).
Note in particular that this result also applies to Gabor Riesz sequences.

\begin{theorem} \label{dualgaborlocalization}
Let $g \in M^1$ and $\Lambda \subset \R^{2d}$ be such that
$\Gc(g,\Lambda)$ is a Gabor frame sequence in $L^2(\R^d)$,
with canonical dual frame sequence
$\tGc = \set{\tg_\lambda}_{\lambda \in \Lambda}$.
Then the following statements hold:

\begin{enumerate}
\item[(a)]
$\tg_\lambda \in M^1$ for all $\lambda \in \Lambda$,

\medskip
\item[(b)] $\sup_\lambda \norm{\tg_\lambda}_{M^1} < \infty$, and

\medskip
\item[(c)]
$\tGc$ is a set of Gabor molecules with respect to an envelope
$\Gamma \in W(\Cc,\ell^1)$.
\end{enumerate} 

\smallskip\noindent
Furthermore, the same conclusions hold when $\tGc$ is replaced by
the canonical Parseval frame $S^{-1/2}(\Gc(g,\Lambda))$.
\end{theorem}

\begin{proof}
(a) Since $g \in M^1$, we have by Theorem~\ref{M1local}
that $(\Gc(g,\Lambda),a)$ is $\ell^1$-self-localized.
Theorem~\ref{selflocthm} therefore implies that $(\tGc,a)$ is
$\ell^1$-self-localized as well.
Hence, by definition, there exists $r \in \ell^1(G)$ such that
$$|\ip{\tg_\lambda}{\tg_\mu}|
\LT r_{a(\lambda) - a(\mu)}.$$
Since $\Gc(g,\Lambda)$ is a Bessel sequence, $\Lambda$ has finite
density by Theorem~\ref{gaborfinite}.
Thus $K = \sup_\lambda |a^{-1}(\lambda)| < \infty$.

For each $\lambda \in \Lambda$, the frame expansion of $\tg_\lambda$ is
\begin{equation} \label{glambdaexpand}
\tg_\lambda
\EQ \sum_{\mu \in \Lambda} \ip{\tg_\lambda}{\tg_\mu} \, g_\mu,
\end{equation}
with convergence in $L^2(\R^d)$.
However,
$\set{\ip{\tg_\lambda}{\tg_\mu}}_{\mu \in \Lambda} \in \ell^1(\Lambda)$,
so by \cite[Thm.~12.1.8]{Gro01} the series on the right-hand side
of \eqref{glambdaexpand} converges in $M^1$-norm,
and therefore $\tg_\lambda \in M^1$.

\medskip
(b) Since translation and modulation are isometries in $M^1$-norm,
we have
\begin{align*}
\norm{\tg_\lambda}_{M^1}
& \LE \sum_{\mu \in \Lambda} |\ip{\tg_\lambda}{\tg_\mu}| \, \norm{g_\mu}_{M^1}
      \\
& \EQ \sum_{j \in G} \sum_{\mu \in a^{-1}(j)}
      |\ip{\tg_\lambda}{\tg_\mu}| \, \norm{g}_{M^1}
\LE K \, \norm{r}_{\ell^1} \, \norm{g}_{M^1}.
\end{align*}

\medskip
(c) Set
\begin{align*}
\Q_{\alpha,\beta}(x,y)
& \EQ [0,\alpha)^d \times [0,\beta)^d + (x,y), \\[1 \jot]
\R_{\alpha,\beta}(x,y)
& \EQ [-\alpha,\alpha)^d \times [-\beta,\beta)^d + (x,y),
\end{align*}
and define
$$\Gamma(z)\
\EQ K \sum_{j \in G} \, r_j \sup_{w \in \R_{\alpha,\beta}(j)}
    |V_\gamma g(z-w)|,
\qquad z \in \R^{2d}.$$
Because $r \in \ell^1$ and $V_\gamma g \in W(\Cc,\ell^1)$, we have
$\Gamma \in W(\Cc,\ell^1)$ as well.

Given $\mu \in \Lambda$, recall that $a(\mu) = \bar\mu$, so
$\mu = a(\mu) + \tilde\mu \in \Q_{\alpha\beta}(a(\mu))$.
Also,
$$\Q_{\alpha\beta}(j + a(\lambda))
\EQ [0,\alpha)^d \times [0,\beta)^d + j + \lambda - \tilde\lambda
\SUBSET \R_{\alpha\beta}(j + \lambda).$$
Therefore, taking the STFT of both sides of \eqref{glambdaexpand}, we have
\begin{align*}
|V_\gamma \tg_\lambda(z)|
& \LE \sum_{\mu \in \Lambda} |\ip{\tg_\lambda}{\tg_\mu}| \, |V_\gamma g_\mu(z)|
      \\[1 \jot]
& \LE \sum_{j \in G} \sum_{\mu \in a^{-1}(j)} r_{a(\mu) - a(\lambda)} \,
      |V_\gamma g(z-\mu)|
      \allowdisplaybreaks \\[1 \jot]
& \LE K \, \sum_{j \in G} r_{j - a(\lambda)} \,
      \sup_{w \in \Q_{\alpha\beta}(j)} |V_\gamma g(z-w)|
      \allowdisplaybreaks \\[1 \jot]
& \EQ K \, \sum_{j \in G} r_j \,
      \sup_{w \in \Q_{\alpha\beta}(j + a(\lambda))} |V_\gamma g(z-w)|
      \allowdisplaybreaks \\[1 \jot]
& \LE K \, \sum_{j \in G} r_j \,
      \sup_{w \in \R_{\alpha\beta}(j + \lambda)} |V_\gamma g(z-w)|
\EQ \Gamma(z-\lambda).
\end{align*}
Thus $\tGc$ is a set of Gabor molecules.

Finally, by Theorem~\ref{selflocthm}, the canonical Parseval frame
$S^{-1/2}(\Gc(g,\Lambda))$ is $\ell^1$-self-localized,
and computations similar to the ones above extend the results to
the Parseval frame.
\end{proof}

\smallskip
\subsection{Gabor Molecules} \label{molecules}

We close by noting that many of the results of the preceding sections
proved for Gabor frames carry over to frames of Gabor molecules with
only minor changes in proof.

\begin{theorem}
Let $\Fc = \set{f_\lambda}_{\lambda \in \Lambda}$ be a set of Gabor molecules
with respect to an envelope $\Gamma \in W(\C,\ell^2)$.
Let $\phi \in M^1$ and $\alpha$, $\beta > 0$ be given.
Then the following statements hold.

\smallskip
\begin{enumerate}
\item[(a)]
If $1 \le p \le 2$ and $\Gamma \in W(\Cc,\ell^p)$,
then $(\Fc, a, \Gc(\phi, G))$ is $\ell^p$-localized.
Further, $f_\lambda \in M^p$ for every $\lambda \in \Lambda$.

\medskip
\item[(b)]
If $\Gamma \in W(\Cc,\ell^1)$, then
$(\Fc,a)$ is $\ell^1$-self-localized.
\end{enumerate}
\end{theorem}

\begin{theorem}
Let $\Fc = \set{f_\lambda}_{\lambda \in \Lambda}$
be a set of Gabor molecules with respect to an envelope
$\Gamma \in W(\Cc,\ell^2)$.
If $\Fc$ is a frame for $L^2(\R^d)$ then the following statements hold.

\smallskip
\begin{enumerate}
\item[(a)]
$1 \le \BD^-(\Lambda) \le \BD^+(\Lambda) < \infty$.

\medskip
\item[(b)]
$\cM^-(\Fc) = \frac1{\BD^+(\Lambda)}$
and
$\cM^+(\Fc) = \frac1{\BD^-(\Lambda)}$.

\medskip
\item[(c)]
If $\BD^+(\Lambda) > 1$ then $\Fc$ has infinite excess, and
there exists an infinite subset $J \subset \Lambda$ such that
$\set{f_\lambda}_{\lambda \in \Lambda \setminus J}$
is a frame for $L^2(\R^d)$.

\medskip
\item[(d)]
If $\Gamma \in W(\Cc,\ell^1)$ and $\BD^-(\Lambda) > 1$, then there exists
$J \subset \Lambda$ with $\BD^+(J) = \BD^-(J) > 0$ such that
$\set{f_\lambda}_{\lambda \in \Lambda \setminus J}$
is a frame for $L^2(\R^d)$.

\medskip
\item[(e)] 
If $\Gamma\in W(\Cc,\ell^1)$ then the canonical dual frame
$\tilde{\Fc}= \set{\tf_\lambda}_{\lambda \in \Lambda}$
is a set of Gabor molecules with respect to an envelope function
$\tilde{\Gamma} \in W(\Cc,\ell^1)$.
\end{enumerate}
\end{theorem}

\smallskip
\section{Relations Among the Localization and Approximation Properties}
\label{relationsec}

\begin{figure}[ht]
\scalebox{.6}{\includegraphics{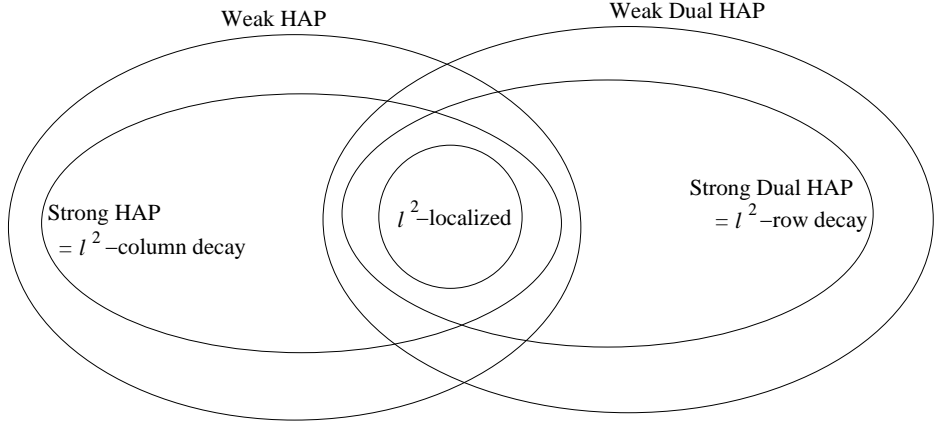}}
\caption{Relations among the localization and approximation properties
for $p=2$, under the assumptions that $\Fc$, $\Ec$ are frames and
$D^+(I,a)<\infty$. \label{fig1}}
\end{figure}

We conclude this work by determining the relationships that hold among the
localization and approximation properties described in section \ref{section2}.
For the case that $\Fc$ and $\Ec$ are both frames for $H$ and 
the upper density $D^+(I,a)$ is finite,
these relations can be summarized in the diagram in Figure~\ref{fig1}.

\smallskip
\subsection{Implications Among the Localization and Approximation Properties}

\begin{theorem} \label{relations}
Let $\Fc = \set{f_i}_{i \in I}$ and
$\Ec = \set{e_j}_{j \in G}$ be sequences in $H$, and
let $a : I \to G$ be an associated map.
Then the following statements hold.

\smallskip
\begin{enumerate}
\item[(a)] If $\Fc$ is a frame for $H$, then
$\ell^2$-column decay implies the strong HAP.

\medskip
\item[(b)] If $\Fc$ is a frame for $H$ and
$\sup_j \norm{e_j} < \infty$, then
the strong HAP implies $\ell^2$-column decay.

\medskip
\item[(c)] If $\Ec$ is a frame for $H$, then
$\ell^2$-row decay implies the strong dual HAP.

\medskip
\item[(d)] If $\Ec$ is a frame for $H$ and
$\sup_i \norm{f_i} < \infty$, then
the strong dual HAP implies $\ell^2$-row decay.

\medskip
\item[(e)] If $\Fc$ is a frame for $H$, then
the strong HAP implies the weak HAP.
If $\Fc$ is a Riesz basis for $H$, then
the weak HAP implies the strong HAP.

\medskip
\item[(f)] If $\Ec$ is a frame for $H$, then
the strong dual HAP implies the weak dual HAP.
If $\Ec$ is a Riesz basis for $H$, then
the weak dual HAP implies the strong dual HAP.

\medskip
\item[(g)] If $D^+(I,a) < \infty$ and $1 \le p < \infty$, then
$\ell^p$-localization implies both $\ell^p$-column and $\ell^p$-row decay.
\end{enumerate}
\end{theorem}

Before giving the proof of Theorem~\ref{relations}, we note that
in Section~\ref{relationappend} we construct counterexamples to most of
the converse implications of Theorem~\ref{relations}, including the following.

\begin{enumerate}
\item[(a)]
There exist orthonormal bases $\Ec$, $\Fc$ such that
$(\Fc,a,\Ec)$ does not have $\ell^2$-column decay, and hence does
not satisfy the strong HAP.

\medskip
\item[(b)]
There exists a frame $\Fc$ and orthonormal basis $\Ec$ such that
$(\Fc,a,\Ec)$ satisfies the weak HAP but not the strong HAP.

\medskip
\item[(c)]
There exists a frame $\Fc$ and orthonormal basis $\Ec$ such that
$D^+(I,a) < \infty$, $(\Fc,a,\Ec)$ has both $\ell^2$-column decay and
$\ell^2$-row decay, but fails to have $\ell^2$-localization.

\medskip
\item[(d)]
There exists a Riesz basis $\Fc$ and orthonormal basis $\Ec$ such that
$(\Fc,a,\Ec)$ has $\ell^2$-column decay but not $\ell^2$-row decay.
\end{enumerate}

\begin{proof}[Proof of Theorem~\ref{relations}]
(a) Assume $\Fc$ is a frame for $H$, with frame bounds $A$, $B$, and
suppose that $(\Fc,a,\Ec)$ has $\ell^2$-column decay.
For each $j \in G$ we have $e_j = \sum_{i \in I} \ip{e_j}{f_i} \, \tf_i$, so
$$\biggnorm{e_j - \sum_{i \in I_{N_\eps(j)}} \ip{e_j}{f_i} \, \tf_i}^2
\EQ \biggnorm{\sum_{i \in I \setminus I_{N_\eps(j)}} \ip{e_j}{f_i} \, \tf_i}^2
\LE \frac1A \, \sum_{i \in I \setminus I_{N_\eps(j)}} |\ip{e_j}{f_i}|^2
\LT \frac{\eps}A.$$
Consequently, $\ell^2$-column decay implies the strong HAP.

\medskip
(b) Assume that $\Fc$ is a frame for $H$,
that $\Ec$ is uniformly bounded above in norm,
and that $(\Fc,a,\Ec)$ has the strong HAP.
Let $S$ be the frame operator for $\Fc$.
Since $S$ is bounded, $C = \sup_j \norm{Se_j} < \infty$.
Let $\tFc = \set{\tf_i}_{i \in I}$ be the canonical dual frame to $\Fc$.
Since $\tf_i = S^{-1} f_i$ and $S$ is self-adjoint, we have
\begin{align*}
\sum_{i \in I \setminus I_{N_\eps(j)}} |\ip{e_j}{f_i}|^2
& \EQ \Bigip{e_j}
      {\sum_{i \in I \setminus I_{N_\eps(j)}} \ip{e_j}{f_i} \, f_i}
      \\
& \EQ \Bigip{Se_j}
            {\sum_{i \in I \setminus I_{N_\eps(j)}} \ip{e_j}{f_i} \, \tf_i}
      \allowdisplaybreaks \\
& \LE \norm{Se_j} \,
      \biggnorm{\sum_{i \in I \setminus I_{N_\eps(j)}} \ip{e_j}{f_i} \, \tf_i}
      \\
& \LE C \, \biggnorm{e_j - \sum_{i \in I_{N_\eps(j)}} \ip{e_j}{f_i} \, \tf_i}
\LT C \eps.
\end{align*}
Consequently, $(\Fc,a,\Ec)$ has $\ell^2$-column decay.

\medskip
(c), (d) These arguments are entirely symmetrical to the ones for (a), (b).

\medskip
(e) Clearly the strong HAP trivially implies the weak HAP.

Suppose that $\Fc$ is a Riesz basis for $H$, and that $(\Fc,a,\Ec)$
satisfies the weak HAP.
Since $\tFc$ is also a Riesz basis for $H$,
there exist constants $A'$, $B'$ such that
$$A' \, \sum_{i \in I} |a_i|^2
\LE \biggnorm{\sum_{i \in I} a_i \tf_i}^2
\LE B' \, \sum_{i \in I} |a_i|^2$$
for any square-summable sequence of scalars $(a_i)$.
Fix any $\eps > 0$, and let $c_{i,j}$ be the
numbers from \eqref{weakHAPdef}.
Then for any $j \in G$,
\begin{align*}
\eps
\GT \biggnorm{e_j - \sum_{i \in N_\eps(j)} c_{i,j} \tf_i}^2
& \EQ \biggnorm{e_j - \sum_{i \in N_\eps(j)} \ip{e_j}{f_i} \, \tf_i +
      \sum_{i \in N_\eps(j)} \bigparen{\ip{e_j}{f_i} - c_{i,j}} \, \tf_i}^2
      \\[1 \jot]
& \EQ \biggnorm{\sum_{i \in I \setminus N_\eps(j)} \ip{e_j}{f_i} \, \tf_i +
                \sum_{i \in N_\eps(j)} (\ip{e_j}{f_i} - c_{i,j}) \, \tf_i}^2
      \allowdisplaybreaks \\[1 \jot]
& \GE A' \,
      \biggparen{\sum_{i \in I \setminus N_\eps(j)} |\ip{e_j}{f_i}|^2 +
      \sum_{i \in N_\eps(j)} |\ip{e_j}{f_i} - c_{i,j}|^2}
      \allowdisplaybreaks \\[1 \jot]
& \GE A' \, \sum_{i \in I \setminus N_\eps(j)} |\ip{e_j}{f_i}|^2
      \allowdisplaybreaks \\[1 \jot]
& \GE \frac{A'}{B'} \,
      \biggnorm{\sum_{i \in I \setminus N_\eps(j)} \ip{e_j}{f_i} \, \tf_i}^2
      \\[1 \jot]
& \EQ \frac{A'}{B'} \,
      \biggnorm{e_j - \sum_{i \in N_\eps(j)} \ip{e_j}{f_i} \, \tf_i}^2.
\end{align*}
Hence $(\Fc,a,\Ec)$ satisfies the strong HAP.

\medskip
(f) This argument is symmetrical to the one for (e).

\medskip
(g) Assume that $(\Fc,a,\Ec)$ is $\ell^p$-localized and that
$D^+(I,a) < \infty$.
Then we have $K = \sup_{n \in G} |a^{-1}(n)| < \infty$.
By definition of $\ell^p$-localization, there exists an
$r \in \ell^p(G)$ such that $|\ip{f_i}{e_j}| \le r_{a(i) - j}$
for all $i \in I$ and $j \in G$.
Given $\eps>0$,
let $N_\eps$ be such that
$$\sum_{\ell \in G \setminus S_{N_\eps}(0)} r_\ell^p \LT \eps.$$
Then
$$\sum_{i \in I \setminus I_{N_\eps}(j)} |\ip{f_i}{e_j}|^p
\LE \sum_{n \in G \setminus S_{N_\eps}(j)} \,
    \sum_{i \in a^{-1}(n)} r_{a(i) - j}^p
\LE K \sum_{n \in G \setminus S_{N_\eps}(j)} r_{n - j}^p
\LE K \eps.$$
Thus $(\Fc,a,\Ec)$ has $\ell^p$-column decay.
Additionally,
$$\sum_{j \in G \setminus S_{N_\eps}(a(i))} |\ip{f_i}{e_j}|^p
\LE \sum_{j \in G \setminus S_{N_\eps}(a(i))} r_{a(i)-j}^p
\LE \sum_{\ell \in G \setminus S_{N_\eps}(0)} r_\ell^p
\LE \eps,$$
so $(\Fc,a,\Ec)$ has $\ell^p$-row decay as well.
\end{proof}

\smallskip
\subsection{Counterexamples}
\label{relationappend}

In this section, we provide examples showing that most of the
implications of Theorem~\ref{relations} are sharp, along with
several other useful examples.

The following example constructs orthonormal bases $\Fc$ and $\Ec$
such that $(\Fc,a,\Ec)$ does not satisfy the strong HAP.

\begin{example} \label{noHAP}
For each $n \in \N$, let $H_n$ be an $n$-dimensional Hilbert space with
orthonormal basis $\set{e_j^n}_{j=1}^n$.
Let $H = \ell^2 \oplus \sum_{n=1}^\infty H_n$, the orthogonal direct
sum of $\ell^2$ and the $H_n$.
Let $\set{e_j}_{j \le 0}$ be an orthonormal basis for $\ell^2$, and
let $\set{e_j}_{j > 0}$ be the orthonormal bases for the $H_n$ put into
their natural order, i.e., $e_{\frac{(n-1)n}2+j} = e_j^n$ for
$j = 1, \dots, n$.
Then $\Ec = \set{e_j}_{j \in \Z}$ is an orthonormal basis for~$H$.

Let $\omega_n = e^{2\pi i/n}$ be a primitive $n$th root of unity, and
let $\set{f_k^n}_{k=1}^n$ be the harmonic orthonormal basis for $H_n$ given by
$$f_k^n \EQ \frac1{\sqrt{n}} \sum_{j=1}^n \omega_n^{jk} e_j^n,
\qquad k = 1, \dots, n.$$
Let $\set{f_k}_{k > 0}$ be the $f_k^n$ put in their natural order,
and for $k \le 0$ set $f_k=e_k$.
Then $\Fc = \set{f_k}_{k \in \Z}$ is an orthonormal basis for $H$.

Let $a \colon \Z \to \Z$ be the identity map.
Fix any $N > 0$.
If $j = \frac{(n-1)n}2 + 1$, then
$\ip{f_k}{e_j} = 0$ for all $k < j$ or $k \ge j+n$.
Since $I_{2N}(j) = S_{2N}(j) = [j-N,j+N) \cap \Z$,
we therefore have for $n>N$ that
\begin{align*}
\sum_{k \in \Z \setminus I_{2N}(j)} |\ip{f_k}{e_j}|^2
& \EQ \sum_{k=j+N}^{j+n-1} |\ip{f_k}{e_j}|^2 \\
& \EQ \sum_{k=N+1}^n |\ip{f_k^n}{e_1^n}|^2
\EQ \sum_{k=N+1}^n \frac1n
\EQ \frac{n-N-1}n.
\end{align*}
This quantity approaches $1$ as $n \to \infty$, so
$(\Fc,a,\Ec)$ fails to have $\ell^2$-column decay, and
hence by Theorem~\ref{relations} also fails the strong HAP.
~\qed
\end{example}

The following example shows that the weak HAP need not imply the strong HAP
if $\Fc$ is not a Riesz basis
(compare to part~(e) of Theorem~\ref{relations}).
Note that in this example, $\Ec$ is actually an orthonormal basis for $H$.

\begin{example}
Let $\Ec$ and $\Fc$ be as in Example~\ref{noHAP}.
Define
$$g_{2i}^n \EQ \frac1{\sqrt{2}} e_i^n, \qquad
g_{2i-1}^n \EQ \frac1{\sqrt{2}} f_i^n, \qquad
i = 1, \dots, n.$$
Let $\set{g_i}_{i > 0}$ be the $g_i^n$ put in their natural order,
and for $i \le 0$ set $g_i = e_i$.
Then $\Gc = \set{g_i}_{i \in \Z}$ is a Parseval frame for $H$,
i.e., the frame bounds are $A=B=1$.
In particular, $\Gc$ is its own dual frame.

Define $a \colon \Z \to \Z$ by
$a(i) = i$ for $i \le 0$ and
$$a\bigparen{\tfrac{(2n-1)2n}2 + 2i-1}
\EQ a\bigparen{\tfrac{(2n-1)2n}2 + 2i}
\EQ \tfrac{(n-1)n}2 + i,
\qquad i = 1, \dots, 2n,$$
i.e., $a$ associates the elements $g_{2i-1}^n$ and $g_{2i}^n$ of $\Gc$
with the element $e_i^n$ in $\Ec$.

Given any $j \in \Z$, we have
$g_j = e_j$ for $j \le 0$, and
$\sqrt2 g_{2j}^n = e_j^n$ for $j > 0$,
so clearly $(\Gc,a,\Ec)$ satisfies the weak HAP.
However, given any $N>0$ we have
$$\sum_{i \in \Z \setminus I_{2N}(j)} |\ip{g_i}{e_j}|^2
\GE \frac12 \sum_{k \in \Z \setminus S_{2N}(j)} |\ip{f_k}{e_j}|^2,$$
and as we saw in Example~\ref{noHAP}, we cannot make this quantity
arbitrarily small independently of $j$.
Thus $(\Gc,a,\Ec)$ fails the strong HAP.
~\qed
\end{example}

The following example shows that the assumption of $\ell^2$-localization
alone does not guarantee that the upper density is finite.
In particular, this shows that the hypothesis in \cite[Thm.~3.3]{BCHL05a}
that $\inf_i \norm{f_i} > 0$ is necessary.

\begin{example}
Let $\Ec = \set{e_n}_{n \in \Z}$ be an orthonormal basis for $H$
and let
$\Fc = \set{e_n}_{n > 0} \cup \set{2^n e_0}_{n \le 0}$.
Note that $\Fc$ is a frame sequence in $H$.
If we let $a(i) = i$ for $i > 0$ and $a(i) = 0$ for $i \le 0$,
then $|a^{-1}(0)| = \infty$,
so $D^+(I,a) = \infty$.
On the other hand,
$\sup_{i \in \Z} |\ip{f_i}{e_{j+a(i)}}| = 1$ if $j = 0$
and $0$ otherwise, so $(\Fc,a,\Ec)$ is $\ell^2$-localized.
~\qed
\end{example}

The following example shows that the converse of part~(g) of
Theorem~\ref{relations} fails in general, i.e.,
$\ell^2$-column decay combined with $\ell^2$-row decay does not imply
$\ell^2$-localization, even if $D^+(I,a) < \infty$.

\begin{example}
Let $\Ec =\set{e_j}_{j \in \Z}$ be an orthonormal basis for $H$,
and define $\Fc = \set{f_j}_{j \in \Z}$ by
$$f_j \EQ e_j + \biggparen{\frac{1}{4+|j|}}^{1/2} e_{-j},
\qquad j \in \Z$$
Let $a \colon \Z \to \Z$ be the identity map.
Then $D^+(a,\Z) = 1$, and $I_N(j) = S_N(j)$ for all $j$ and $N$.
For $j \ne 0$, we have
$$\sup_{i \in \Z} |\ip{f_i}{e_{2j+i}}|^2
\EQ |\ip{f_{-j}}{e_j}|^2 
\EQ \frac{1}{4+|j|},$$
so $(\Fc,a,\Ec)$ is not $\ell^2$-localized.
On the other hand, since $\ip{f_i}{e_j} \ne 0$ only when $i = \pm j$, we have
\begin{equation} \label{calculation}
\sum_{i \in \Z \setminus S_{N_\eps}(j)} |\ip{f_i}{e_j}|^2
\EQ \begin{cases}
       0, & -\frac{N}4 < j \le \frac{N}4, \\
       \frac1{4+|j|}, & \text{otherwise}.
    \end{cases}
\end{equation}
By taking $N_\eps$ large enough, we can make this quantity arbitrarily
small, independently of~$j$.
Thus $(\Fc,a,\Ec)$ has $\ell^2$-column decay,
and a similar argument shows it has $\ell^2$-row decay.

Note that no other choice for the map $a$ would help in this example,
for if $(\Fc,a,\Ec)$ has $\ell^2$-column decay, then
$\sup_j |a(j)-j| < \infty$.
Thus $a$ can only be a bounded perturbation of the identity.
For such an $a$, there always exists an $N$ sufficiently large so that
for $|j| > N$ an inequality similar to \eqref{calculation} will hold.

Note also in this example that $\Fc$ is a frame.
For, given $f \in H$ we have
$$\sum_{j \in \Z} |\ip{f}{f_j-e_j}|^2
\EQ \sum_{j \in \Z} \biggparen{\frac1{4+|j|}} |\ip{f}{e_{-j}}|^2
\LE \frac14 \sum_{j \in \Z} |\ip{f}{e_j}|^2
\EQ \frac14 \, \norm{f}^2,$$
and therefore, by the triangle inequality,
\begin{align*}
\biggparen{\sum_{j \in \Z} |\ip{f}{f_j}|^2}^{1/2}
& \GE \biggparen{\sum_{j \in \Z} |\ip{f}{e_j}|^2}^{1/2} \minus
      \biggparen{\sum_{j \in \Z} |\ip{f}{f_j-e_j}|^2}^{1/2}
      \\[1 \jot]
& \GE \norm{f} - \frac12 \, \norm{f}
\EQ \frac12 \, \norm{f}.
\end{align*}
Thus $\Fc$ has a lower frame bound of $1/4$, and a similar
calculation shows it has an upper frame bound of $9/4$.
~\qed
\end{example}

The following example shows that $\ell^2$-column decay does not imply
$\ell^2$-row decay.
By interchanging the roles of $\Fc$ and $\Ec$ in this example, we also see
that $\ell^2$-row decay does not imply $\ell^2$-column decay.

\begin{example}
Index an orthonormal basis for $H$ as
$\Ec = \Span\set{e_j^n}_{n \in \N, j=1,\dots,n}$, and set
$H_n = \Span\set{e_j^n}_{j=1,\dots,n}$.
Define
$$f_i^n
\EQ \begin{cases}
    e_1^n, & i=1, \\
    \frac1{2\sqrt{n}} e_1^n + e_i^n, & i = 2, \dots, n,
    \end{cases}$$
and
$$\tf_i^n
\EQ \begin{cases}
    e_1^n - \frac1{2\sqrt{n}} \sum_{j=2}^n e_j^n, & i=1, \\
    e_i^n, & i = 2, \dots, n.
    \end{cases}$$
Clearly $f_i^n$, $\tf_i^n \in H_n$, and a straightforward calculation
shows that $\set{f_i^n}_{i=1}^n$ and $\set{\tf_i^n}_{i=1}^n$
are biorthogonal sequences in $H_n$.
Since $H_n$ is $n$-dimensional, this shows that these are dual Riesz
bases for $H_n$.
Given any scalars $\set{a_i}_{i=1}^n$, we have
\begin{align*}
\biggnorm{\sum_{i=1}^n a_i f_i^n}
& \LE \biggnorm{\frac{\sum_{i=1}^n a_i}{2\sqrt{n}} \, e_1^n} +
      \biggnorm{\sum_{i=1}^n a_i e_i^n}
      \\[1 \jot]
& \LE \frac{\sum_{i=1}^n |a_i|}{2\sqrt{n}} +
      \biggparen{\sum_{i=1}^n |a_i|^2}^{1/2}
\LE \frac32 \biggparen{\sum_{i=1}^n |a_i|^2}^{1/2},
\end{align*}
and similarly
$\bignorm{\sum_{i=1}^n a_i f_i^n}
\ge \frac12 \bigparen{\sum_{i=1}^n |a_i|^2}^{1/2}$.
Thus $\set{f_i^n}_{i=1}^n$ has Riesz bounds $\frac12$, $\frac32$.
Since $H$ is the orthogonal direct sum of the $H_n$ and the Riesz bounds are
independent of $n$, we conclude that
$\Fc = \set{f_i^n}_{n \in \N, i=1,\dots,n}$ and
$\tFc = \set{\tf_i^n}_{n \in \N, i=1,\dots,n}$
are dual Riesz bases for $H$.

Another direct calculation shows that
$$|\ip{f_i^m}{e_j^n}|
\EQ \begin{cases}
    1, & i=j, m=n, \\
    \frac1{2\sqrt{n}}, & m=n, j=1, i=2,\dots,n, \\
    0, & \text{otherwise}.
    \end{cases}$$
Consequently, given any $N$, we have for each $n > N$ and $j=1,\dots,n$ that
$$\sum_{m>N} \sum_{i=1}^m |\ip{f_i^m}{e_j^n}|^2
\EQ \sum_{i=1}^n \ip{f_i^n}{e_j^n}|^2
\EQ \frac1{4n},$$
while for $n \le N$ this sum is zero.
Hence, by taking $N$ large enough this sum is less than $\eps$ independently
of $n \in \N$ and $j=1, \dots, n$.
Thus, with $a$ as the identity map, $(\Fc,a,\Ec)$ has $\ell^2$-row decay.
On the other hand, if $m>N$ then we have for each $i=1,\dots,n$ that
$$\sum_{n>N} \sum_{j=1}^n |\ip{f_i^m}{e_j^n}|^2
\EQ \sum_{j=1}^n |\ip{f_i^n}{e_j^n}|^2
\EQ \sum_{j=2}^n \frac1{4n}
\EQ \frac{n-1}{4n}.$$
After mapping the index set of $\Ec$ and $\Fc$ onto $\Z$,
similarly to Example~\ref{noHAP},
we conclude that $(\Fc,a,\Ec)$ does not have $\ell^2$-row decay.
~\qed
\end{example}

The following example illustrates the importance of the map $a$ in
determining localization properties.

\begin{example} \label{selflocexample}
Let $\Fc = \set{f_n}_{n \in \Z}$ be an orthonormal basis for $H$,
and define $a \colon \Z \to \Z$ by $a(2n) = a(2n+1) = n$.
Then $(\Fc,a)$ is $\ell^1$-self-localized, and by Example~\ref{specialcases}
we have $\cM^\pm(\Fc) = 1$.
However, $D^\pm(I,a) = 2$.
Hence, by \cite[Thm.~3.5(c)]{BCHL05a}, there cannot be any Riesz basis $\Ec$
such that $(\Fc,a,\Ec)$ has both $\ell^2$-column decay and $\ell^2$-row decay.
In particular $(\Fc,a,\Ec)$ cannot be $\ell^2$-localized for any Riesz
basis $\Ec$, and $(\Fc,a,\Fc)$ is not $\ell^2$-localized.

However, if we let $\Ec = \set{f_{2n}}_{n \in \Z}$, then $\Ec$ is a
Riesz sequence (but not a Riesz basis), and
$(\Fc,a,\Ec)$ is $\ell^1$-localized.
Since $\Ec$ is a Riesz sequence and $\clspan(\Fc) = H$, we have
$\cM^\pm(\Ec) = \cM(\Ec;p,c) = 1$ by Example~\ref{specialcases}(b).
On the other hand, since
$P_\Ec(f_{2n}) = 1$ and $P_\Ec(f_{2n+1}) = 0$,
it follows directly that $\cM_\Ec(\Fc;p,c) = \frac12$.
Thus $\cM_\Ec(\Fc;p,c) \, D(p,c) = 1$,
in accordance with \cite[Thm.~3.4]{BCHL05a}.
~\qed
\end{example}

The following example shows that $\ell^1$-localization with respect to
the canonical dual frame does not imply $\ell^1$-self-localization.

\begin{example} \label{selfexample}
Let $\Ec = \set{e_i}_{i \in \Z}$ be an orthonormal basis for $H$,
and let $a \colon \Z \to \Z$ be the identity map.
Fix $\frac12 < c_0 < 1$, and for $i \ne 0$ choose $c_i > 0$ in such a
way that
$$\sum_{i \in \Z} c_i^2 \EQ 1
\qquad\text{and}\qquad
\sum_{i \in \Z} c_i \EQ \infty.$$
Define
$$f_0 \EQ \sum_{i \in \Z} c_i e_i
\qquad\text{and}\qquad
f_i = e_i \text{ for } i \ne 0.$$
If we set $T(e_i) = f_i$, then $T$ extends to a bounded mapping on $H$.
Further, if $f = \sum_i \ip{f}{e_i} \, e_i \in H$, then
$$\norm{(\one - T)f}^2
\EQ \bignorm{\ip{f}{e_0} \, (e_0 - f_0)}^2
\EQ |\ip{f}{e_0}|^2 \, \Bigparen{|1 - c_0|^2 + \sum_{i \ne 0} c_i^2}
\LE (2 - 2c_0) \, \norm{f}^2,$$
so
$\norm{\one - T} \le 2 - 2c_0 < 1$.
Hence $T$ is a continuous bijection of $H$ onto itself, so
$\Fc = \set{f_i}_{i \in \Z}$ is a Riesz basis for $H$.
Therefore $\ip{f_i}{\tf_j} = \delta_{ij}$, so $(\Fc,a)$ is $\ell^1$-localized
with respect to its dual frame.
However, $\ip{f_0}{f_j} = c_j$, so $(\Fc,a)$ is not $\ell^1$-self-localized.
~\qed
\end{example}

\smallskip
\appendix
\section{Ultrafilters} \label{ultraappend}

In this appendix we provide a brief review of ultrafilters and their
basic properties.
For additional information, we refer to \cite[Chap.~3]{HS98}.
Filters were introduced by H.~Cartan \cite{Car37a}, \cite{Car37b}
in order to characterize continuous functions on general topological spaces.
Soon after, it was realized that the set of ultrafilters endowed with the
proper topology is the Stone-\v{C}ech compactification of a discrete
(or more generally, a completely regular) topological space.
In the following we will restrict our attention to ultrafilters over
the natural numbers $\N$.

\begin{definition}
A collection $p$ of subsets of $\N$ is a \emph{filter} if:

\smallskip
\begin{enumerate}
\item[(a)] $\emptyset \notin p$,

\smallskip
\item[(b)] if $A$, $B \in p$ then $A \cap B \in p$,

\smallskip
\item[(c)] if $A \in p$ and $A \subset B \subset \N$, then $B \in p$.
\end{enumerate}
A filter $p$ is an \emph{ultrafilter} if it is maximal in the sense that:

\smallskip
\begin{enumerate}
\item[(d)] if $p'$ is a filter on $\N$ such that $p \subset p'$, then $p'=p$,
\end{enumerate}
or, equivalently, if

\smallskip
\begin{enumerate}
\item[(d')] for any $A \subset \N$, either $A \in p$ or
$\N \setminus A \in p$
(but not both, because of properties~a and~b).
\end{enumerate}
The set of ultrafilters is denoted by $\beta\N$.
~\qed
\end{definition}

\begin{definition}
Given any $n \in \N$,
$e_n =\set{A \subset \N : n \in A}$
is an ultrafilter, called a \emph{principal ultrafilter}.
It is straightforward to show that any ultrafilter $p$ that contains a finite
set must be one of these principal ultrafilters.
An ultrafilter which contains no finite sets is called \emph{free}.
The set of free ultrafilters is denoted by $\N^*$.
~\qed
\end{definition}

Our main use for ultrafilters is that they provide
a notion of convergence for arbitrary sequences.

\begin{definition}
Let $p \in \beta\N$ be an ultrafilter.
Then we say that a sequence $\set{c_k}_{k \in \N}$ of complex numbers
\emph{converges to $c \in \C$ with respect to $p$} if
for every $\eps > 0$ there exists a set $A \in p$ such that 
$|c_k - c| < \eps$ for all $k \in A$.
In this case we write $c = \plim_{k \in \N} c_k$ or simply $c = \plim c_k$.
~\qed
\end{definition}

The following proposition summarizes the basic properties of convergence
with respect to an ultrafilter.

\begin{proposition}
Let $p \in \beta\N$ be an ultrafilter.
Then the following statements hold.

\smallskip
\begin{enumerate}
\item[(a)] Every bounded sequence of complex scalars $\set{c_k}_{k \in \N}$
converges with respect to $p$ to some $c \in \C$.

\medskip
\item[(b)] $p$-limits are unique.

\medskip
\item[(c)] If $p = e_n$ is a principal ultrafilter, then
$\plim c_k = c_n$.

\medskip
\item[(d)] If $\set{c_k}_{k \in \N}$ is a convergent sequence in the
usual sense, $p$ is a free ultrafilter, and
$\lim_{k \to \infty} c_k = c$, then $\plim c_k = c$.

\medskip
\item[(e)]
If $\set{c_k}_{k \in \N}$ is a bounded sequence and $p$ is a free ultrafilter,
then $\plim_{k \in \N} c_k$ is an accumulation point of
$\set{c_k}_{k \in \N}$.

\medskip
\item[(f)]
If $c$ is an accumulation point of a bounded sequence $\set{c_k}_{k \in \N}$,
then there exists a free ultrafilter $p$ such that
$\plim c_k = c$.
In particular, there exists an ultrafilter $p$ such that
$\plim c_k = \limsup c_k$,
and there exists an ultrafilter $q$ such that
$\qlim c_k = \liminf c_k$.

\medskip
\item[(g)]
$p$-limits are linear, i.e.,
$\plim (a c_k + b d_k) = a \plim c_k + b \plim d_k$.

\medskip
\item[(h)]
$p$-limits respect products, i.e.,
$\plim (c_k d_k) = \bigparen{\plim c_k} \, \bigparen{\plim d_k}$.
\end{enumerate}
\end{proposition}

\smallskip
\section*{Acknowledgments}
We gratefully acknowledge conversations with Karlheinz Gr\"ochenig and
Massimo Fornasier on localization of frames, and thank them for providing
us with preprints of their articles.
We thank
Hans Feichtinger,
Norbert Kaiblinger,
Gitta Kutyniok,
and
Henry Landau
for conversations regarding the details of our arguments,
and also acknowledge helpful conversations with
Akram Aldroubi,
Carlos Cabrelli,
Mark Lammers,
Ursula Molter,
and
Kasso Okoudjou.


\begin{thebibliography}{BCHL03b}

\bibitem[Bag90]{Bag90}
L.~Baggett,
Processing a radar signal and representations of the discrete Heisenberg group,
\textsl{Colloq.\ Math.}, \textbf{60/61} (1990), 195--203.




\bibitem[BCHL03a]{BCHL03a}
R.~Balan, P.~G.~Casazza, C.~Heil, and Z.~Landau,
Deficits and excesses of frames,
\textsl{Adv.\ Comput.\ Math.}, \textbf{18} (2003), 93--116.

\bibitem[BCHL03b]{BCHL03b}
R.~Balan, P.~G.~Casazza, C.~Heil, and Z.~Landau,
Excesses of Gabor frames,
\textsl{Appl.\ Comput.\ Harmon.\ Anal.}, \textbf{14} (2003), 87--106.

\bibitem[BCHL05a]{BCHL05a}
R.~Balan, P.~G.~Casazza, C.~Heil, and Z.~Landau,
Density, overcompleteness, and localization of frames, I. Theory,
preprint (2005).

\bibitem[BCHL05b]{BCHL05b}
R.~Balan, P.~G.~Casazza, C.~Heil, and Z.~Landau,
Density, overcompleteness, and localization of frames,
research announcement (2005).




\bibitem[BCGP02]{BCGP02}
J.~J.~Benedetto, W.~Czaja, and A.~Ya.~Maltsev,
The Balian--Low theorem for the symplectic form on $\R^{2d}$,
\textsl{J. Geom.\ Anal.}, \textbf{13} (2003), 2, 217--232.

\bibitem[BHW95]{BHW95}
J.~J.~Benedetto, C.~Heil, and D.~F.~Walnut,
Differentiation and the Balian--Low theorem,
\textsl{J. Fourier Anal.\ Appl.}, \textbf{1} (1995), 355--402.

\bibitem[BR03]{BR03}
M.~Bownik and Z.~Rzeszotnik,
The spectral function of shift-invariant spaces,
\textsl{Michigan Math.\ J.}, \textbf{51} (2003), 387--414.


\bibitem[Car37a]{Car37a}
H.~Cartan,
Th\'{e}orie des filtres,
\textsl{C. R. Acad.\ Sci.\ Paris}, \textbf{205} (1937), 595--598.

\bibitem[Car37b]{Car37b}
H.~Cartan,
Filtres et ultrafiltres,
\textsl{C. R. Acad.\ Sci.\ Paris}, \textbf{205} (1937), 777--779.

\bibitem[Cas00]{Cas00}
P.~G.~Casazza,
The art of frame theory,
\textsl{Taiwanese J. Math.}, \textbf{4} (2000), 129--201.


\bibitem[CK02]{CK02}
P.~G.~Casazza and N.~J.~Kalton,
Roots of complex polynomials and Weyl--Heisenberg frame sets,
\textsl{Proc.\ Amer.\ Math.\ Soc.}, \textbf{130} (2002), 2313--2318.


\bibitem[Chr03]{Chr03}
O.~Christensen,
``An Introduction to Frames and Riesz Bases,''
Birkh\"auser, Boston, 2003.

\bibitem[CDH99]{CDH99}
O.~Christensen, B.~Deng, and C.~Heil,
Density of Gabor frames,
\textsl{Appl.\ Comput.\ Harmon.\ Anal.}, \textbf{7} (1999), 292--304.

\bibitem[CFZ01]{CFZ01}
O.~Christensen, S.~Favier, and Z.~Felipe,
Irregular wavelet frames and Gabor frames,
\textsl{Approx.\ Theory Appl.\ (N.S.)}, \textbf{17} (2001), 90--101.


\bibitem[Dau90]{Dau90}
I.~Daubechies,
The wavelet transform, time-frequency localization and signal analysis,
\textsl{IEEE Trans.\ Inform.\ Theory}, \textbf{39} (1990), 961--1005.

\bibitem[Dau92]{Dau92}
I. Daubechies,
``Ten Lectures on Wavelets,'' SIAM, Philadelphia, 1992.

\bibitem[DGM86]{DGM86}
I.~Daubechies, A.~Grossmann, and Y.~Meyer,
Painless nonorthogonal expansions,
\textsl{J. Math.\ Phys.}, \textbf{27} (1986), 1271--1283.

\bibitem[DLL95]{DLL95}
I.~Daubechies, H.~Landau, and Z.~Landau,
Gabor time-frequency lattices and the Wexler-Raz identity,
\textsl{J. Fourier Anal.\ Appl.}, \textbf{1} (1995), 437--478.

\bibitem[DH00]{DH00}
B.~Deng and C.~Heil,
Density of Gabor Schauder bases,
in: ``Wavelet Applications in Signal and Image Processing VIII''
(San Diego, CA, 2000), A.~Aldroubi at al., eds.,
Proc.\ SPIE Vol.\ 4119, SPIE, Bellingham, WA, 2000, 153--164.

\bibitem[Fei80]{Fei80}
H.~G.~Feichtinger,
\emph{Banach convolution algebras of {W}iener type},
in: Functions, Series, Operators, Proc.\ Conf.\ Budapest \textbf{38},
Colloq.\ Math.\ Soc.\ J\'anos Bolyai, 1980, 509--524.

\bibitem[Fei81]{Fei81}
H.~G.~Feichtinger,
On a new {S}egal algebra,
\textsl{Monatsh.\ Math.}, \textbf{92} (1981), 269--289.

\bibitem[Fei87]{Fei87}
H.~G.~Feichtinger,
\emph{Banach spaces of distributions defined by decomposition methods,
\textup{II}}, Math.\ Nachr., \textbf{132} (1987), 207--237.

\bibitem[Fei89]{Fei89}
H.~G.~Feichtinger,
Atomic characterizations of modulation spaces through
{G}abor-type representations,
\textsl{Rocky Mountain J. Math.}, \textbf{19} (1989), 113--125.

\bibitem[Fei90]{Fei90}
H.~G.~Feichtinger,
\emph{Generalized amalgams, with applications to {F}ourier transform},
Canad.\ J. Math, \textbf{42} (1990), pp.~395--409.

\bibitem[Fei03]{Fei03}
H.~G.~Feichtinger,
Modulation spaces of locally compact Abelian groups,
in: ``Wavelets and their Applications'' (Chennai, January 2002),
M.~Krishna, R.~Radha and S.~Thangavelu, eds.,
Allied Publishers, New Delhi (2003), 1--56.

\bibitem[Fei92]{Fei92}
H.~G.~Feichtinger,
\emph{Wiener amalgams over {E}uclidean spaces and some of their applications},
in: Function spaces (Edwardsville, IL, 1990),
Lecture Notes in Pure and Appl.\ Math.\ \textbf{136},
Dekker, New York, 1992, 123--137.

\bibitem[FG85]{FG85}
H.~G.~Feichtinger and P.~Gr\"obner,
\emph{Banach spaces of distributions defined by decomposition methods,
\textup{I}}, Math.\ Nachr., \textbf{123} (1985), 97--120.

\bibitem[FG89a]{FG89a}
H.~G.~Feichtinger and K.~Gr\"ochenig,
Banach spaces related to integrable group representations and their
atomic decompositions,~I,
\textsl{J. Funct.\ Anal.}, \textbf{86} (1989), 307--340.

\bibitem[FG89b]{FG89b}
H.~G.~Feichtinger and K.~Gr\"ochenig,
Banach spaces related to integrable group representations and their atomic
decompositions,~II, \textsl{Monatsh.\ Math.}, \textbf{108} (1989), 129--148.

\bibitem[FG97]{FG97}
H.~G.~Feichtinger and K.~Gr\"ochenig,
Gabor frames and time-frequency analysis of distributions,
\textsl{J. Funct.\ Anal.}, \textbf{146} (1997), 464--495.

\bibitem[FZ98]{FZ98}
H.~G.~Feichtinger and G.~Zimmermann,
A Banach space of test functions for Gabor analysis,
in: ``Gabor Analysis and Algorithms: Theory and Applications,''
H.~G.~Feichtinger and T.~Strohmer, eds., Birkh\"auser (1998), 123--170.

\bibitem[For03]{For03}
M.~Fornasier,
Constructive methods for numerical applications in signal processing
and homogenization problems, Ph.D.\ Thesis, U.~Padua, 2003.

\bibitem[Gr\"o93]{Gro93}
K.~Gr\"ochenig,
Irregular sampling of wavelet and short-time Fourier transforms,
\textsl{Constr.\ Approx.}, \textbf{9} (1993), 283--297.

\bibitem[Gr\"o01]{Gro01}
K.~Gr\"ochenig,
``Foundations of Time-Frequency Analysis,'' Birkh\"auser, Boston, 2001.

\bibitem[Gr\"o03]{Gro03}
K.~Gr\"ochenig,
Localized frames are finite unions of Riesz sequences,
\textsl{Adv.\ Comput.\ Math.}, \textbf{18} (2003), 149--157.

\bibitem[Gr\"o04]{Gro04}
K.~Gr\"ochenig,
Localization of frames, Banach frames, and the invertibility of the frame
operator, \textsl{J. Fourier Anal.\ Appl.}, \textbf{10} (2004), 105--132.

\bibitem[GF04]{GF04}
K.~Gr\"ochenig and M.~Fornasier,
Intrinsic localization of frames,
Constructive Approx., to appear (preprint 2004).

\bibitem[GHHK02]{GHHK02}
K.~Gr\"ochenig, D.~Han, C.~Heil, and G.~Kutyniok,
The Balian--Low Theorem for symplectic lattices in higher dimensions,
Appl.\ Comput.\ Harmon.\ Anal., 13 (2002), 169--176.

\bibitem[GL04]{GL04}
K.~Gr\"ochenig and M.~Leinert,
Wiener's lemma for twisted convolution and Gabor frames,
\textsl{J. Amer.\ Math.\ Soc.}, \textbf{17} (2004), 1--18.

\bibitem[GR96]{GR96}
K.~Gr\"ochenig and H.~Razafinjatovo,
On Landau's necessary density conditions for sampling and
interpolation of band-limited functions,
\textsl{J. London Math.\ Soc.\ (2)}, \textbf{54} (1996), 557--565.


\bibitem[HW01]{HW01}
D.~Han and Y.~Wang,
Lattice tiling and the Weyl--Heisenberg frames,
\textsl{Geom.\ Funct.\ Anal.}, \textbf{11} (2001), 742--758.

\bibitem[Hei03]{Hei03}
C.~Heil,
An introduction to weighted Wiener amalgams,
in: ``Wavelets and their Applications'' (Chennai, January 2002),
M.~Krishna, R.~Radha and S.~Thangavelu, eds.,
Allied Publishers, New Delhi (2003), 183--216.

\bibitem[HW89]{HW89}
C.~E.~Heil and D.~F.~Walnut,
Continuous and discrete wavelet transforms,
\textsl{SIAM Review}, \textbf{31} (1989), 628--666.


\bibitem[HS98]{HS98}
N.~Hindman and D.~Strauss,
Algebra in the Stone-\v{C}ech Compactification,
de Gruyter Expositions in Mathematics Vol.~27,
Walter de Gruyter and Co., Berlin, 1998.


\bibitem[Jan94]{Jan94}
A.~J.~E.~M.~Janssen,
Signal analytic proofs of two basic results on lattice expansions,
\textsl{Appl.\ Comput.\ Harmon.\ Anal.}, \textbf{1} (1994), 350--354.

\bibitem[Jan98]{Jan98}
A.~J.~E.~M.~Janssen,
A density theorem for time-continuous filter banks,
in: Signal and image representation in combined spaces,
Y.~Y.~Zeevi and R.~R.~Coifman, eds., Wavelet Anal. Appl., Vol.~7,
Academic Press, San Diego, CA, 1998, 513--523.

\bibitem[Jan03]{Jan03}
A.~J.~E.~M.~Janssen,
On generating tight Gabor frames at critical density,
\textsl{J. Fourier Anal.\ Appl.}, \textbf{9} (2003), 175--214.

\bibitem[JS02]{JS02}
A.~J.~E.~M.~Janssen and T.~Strohmer,
Hyperbolic secants yield Gabor frames,
\textsl{Appl.\ Comput.\ Harmon.\ Anal.}, \textbf{12} (2002), 259--267.



\bibitem[Lan93]{Lan93}
H.~Landau,
On the density of phase-space expansions,
\textsl{IEEE Trans.\ Inform.\ Theory}, \textbf{39} (1993), 1152--1156.


\bibitem[LW03]{LW03}
Y.~Liu and Y.~Wang,
The uniformity of non-uniform Gabor bases,
\textsl{Adv.\ Comput.\ Math.}, \textbf{18} (2003), 345--355.

\bibitem[Lyu92]{Lyu92}
Yu.~I.~Lyubarski\u{i},
 Frames in the Bargmann space of entire functions,
in: ``Entire and subharmonic functions,''
Amer.\ Math.\ Soc., Providence, RI, 1992, 167--180.

\bibitem[RS95]{RS95}
J. Ramanathan and T. Steger,
Incompleteness of sparse coherent states,
\textsl{Appl.\ Comput.\ Harmon.\ Anal.}, \textbf{2} (1995), 148--153.

\bibitem[Rie81]{Rie81}
M. Rieffel,
Von Neumann algebras associated with pairs of lattices in Lie groups,
\textsl{Math.\ Ann.}, \textbf{257} (1981), 403--418.


\bibitem[Sei92]{Sei92}
K.~Seip,
Density theorems for sampling and interpolation in the
Bargmann--Fock space I,
\textsl{J. Reine Angew.\ Math.}, \textbf{429} (1992), 91--106.

\bibitem[SW92]{SW92}
K.~Seip and R.~Wallst\'en,
Sampling and interpolation in the Bargmann--Fock space II,
\textsl{J. Reine Angew.\ Math.}, \textbf{429} (1992), 107--113.


\bibitem[SZ02]{SZ02}
W.~Sun and X.~Zhou,
Irregular wavelet/Gabor frames,
\textsl{Appl.\ Comput.\ Harmon.\ Anal.}, \textbf{13} (2002), 63--76.

\bibitem[SZ03]{SZ03}
W.~Sun and X.~Zhou,
Irregular Gabor frames and their stability,
\textsl{Proc.\ Amer.\ Math.\ Soc.}, \textbf{131} (2003), 2883--2893.


\bibitem[You01]{You01}
R.~Young,
``An Introduction to Nonharmonic Fourier Series,'' Revised First Edition,
Academic Press, San Diego, 2001.

\end{thebibliography}
\end{document}